\documentclass[10pt,reqno]{amsart}%
\usepackage{mathrsfs}
\usepackage{amsfonts,amssymb}
\usepackage{amsmath}
\usepackage{graphicx}
\usepackage{color}
\usepackage[hidelinks]{hyperref}
\newtheorem{thm}{Theorem}[section]
\newtheorem{cor}[thm]{Corollary}
\newtheorem{lem}[thm]{Lemma}
\newtheorem{prop}[thm]{Proposition}

\newtheorem{rem}[thm]{Remark}
\newtheorem{exm}[thm]{Example}

\newtheorem{algo}[thm]{Algorithm}
\numberwithin{equation}{section}

\newcommand\figcaption{\def\@captype{figure}\caption}
\newcommand\tabcaption{\def\@captype{table}\caption}
\newcommand{\norm}[1]{\left\Vert#1\right\Vert}
\newcommand{\abs}[1]{\left\vert#1\right\vert}
\newcommand{\set}[1]{\left\{#1\right\}}
\newcommand{\Real}{\mathbb R}
\newcommand{\mean}[1]{\mathbb{E}\lbrack #1\rbrack}

\newcommand{\meanq}[1]{\mathbb{E}_{Q}\lbrack #1\rbrack}
\newcommand{\meanp}[1]{\mathbb{E}_{P}\lbrack #1\rbrack}
\newcommand{\floor}[1]{\lfloor #1 \rfloor  }
\newcommand{\sfrac}[2]{#1\slash #2}

\newcommand{\cC}{\mathcal{C}}
\newcommand{\cD}{\mathcal{D}}
\newcommand{\cH}{\mathcal{H}}
\newcommand{\cF}{\mathcal{F}}

\newcommand{\cL}{\mathcal{L}}
\newcommand{\cM}{\mathcal{M}}

\newcommand{\cJ}{\mathcal{J}}
\newcommand{\cO}{\mathcal{O}}

\newcommand{\cT}{\mathcal{T}}

\begin{document}
 
\title[Comparison of Wiener chaos and stochastic collocation]{Wiener chaos
vs stochastic collocation methods for linear  advection-diffusion-reaction equations
with multiplicative white noise $^\dagger$}
 \thanks{$^\dagger$ To cite this paper, use \\  
 Z. Zhang, M. V. Tretyakov, B. Rozovskii, and G. E. Karniadakis. Wiener chaos vs stochastic collocation methods for linear advection-diffusion equations with multiplicative white noise. SIAM J. Numer. Anal., 53(1): 153-183, 2015.}

\author{Zhongqiang Zhang}
\thanks{Department of Mathematical Sciences, Worcester Polytechnic Institute,  Worcester, MA, 01609, USA. Email: zzhang7@wpi.edu.}

\author{Michael V. Tretyakov} 
\thanks{School of Mathematical Sciences, University of Nottingham, 
Nottingham, NG7 2RD, UK. 
Email: Michael.Tretyakov@nottingham.ac.uk}

\author{Boris Rozovskii}
\thanks{Division of Applied Mathematics, Brown University, Providence
RI, 02912, USA. Email: boris\_rozovsky@brown.edu.} 

\author{George E. Karniadakis}
\thanks{Division  of Applied Mathematics, Brown University, Providence
RI, 02912, USA. Email: george\_karniadakis@brown.edu.} 

\maketitle

\begin{abstract}
We  compare  Wiener chaos  and stochastic collocation methods   for linear
advection-reaction-diffusion equations
with multiplicative white noise. Both methods are constructed based on a
recursive multi-stage algorithm for long-time integration.
We derive error estimates for both methods
and  compare their numerical performance.
Numerical results confirm that the recursive multi-stage stochastic collocation method
is of order $\Delta$ (time step size) in
the second-order moments while the recursive multi-stage Wiener chaos method
is of order $\Delta^{\mathsf{N}}+\Delta^2$   ($\mathsf{N}$ is the order
of Wiener chaos) for advection-diffusion-reaction equations with commutative noises, in agreement with
the theoretical error estimates.
However, for non-commutative noises, both methods are of order one in the second-order moments.
\end{abstract}

{\em Key words.}
Wong-Zakai approximation,  spectral expansion, multi-stage, weak convergence 

\emph{Mathematics subject classification.} 
Primary 60H15; Secondary 35R60, 60H40

\allowdisplaybreaks[2]
\vskip 10pt

 \textbf{Notation.}

 $q$:  number of Brownian motions (noises).

 $\mathsf{N}$: highest order of Hermite polynomial chaos.

 $\mathsf{n}$: number of basis modes in approximating the Brownian motion.

 $\mathsf{L}$: level of Smolyak sparse grid collocation.

 $\mathsf{M}$: number of Fourier collocation nodes in physical space.

 $\Delta$:  element size (in time) for multi-element spectral approximation of Brownian   motion.

 $\mathsf{K}$: number of elements in time, which is  $\sfrac{T}{\Delta}$ with $T$ the final integration time.

 $\delta t$: time step size for time discretization in the time interval $(0,\Delta]$.

 $\eta(\mathsf{L}, \mathsf{n}q)$: number of sparse grid points at level $\mathsf{L}$ with dimension
 $\mathsf{n}q$.

 \section{Introduction}
Partial differential equations (PDEs) driven by white noise
have different interpretations of stochastic products and lead to different numerical approximations,
unlike the PDEs driven by colored noise.
Specifically, stochastic products for white noise  are usually interpreted with two different products:
the  Ito  product and
the Stratonovich product, see e.g. \cite{Arnold-B74}.
 Though a problem can be  equivalently formulated  using these two products, the use of different products leads to
different performance of numerical solvers for PDEs driven by white noise, especially when  Wiener chaos expansion (WCE)
and stochastic collocation methods (SCM) in random space are used.
In this paper, we will show theoretically and through numerical
examples that for white noise driven PDEs, WCE and SCM have quite different
performance when the noises are commutative.
This is different from how WCE and SCM behave for PDEs driven by colored noise.
For elliptic equations with colored noise, it is demonstrated  in   \cite{BacNTT11,ElmMPT11}
that  there are only small differences in the numerical performance
 of generalized polynomial
chaos expansion  and SCM.

To apply WCE and SCM, we first discretize the Brownian motion with its truncated spectral expansion, see e.g.
\cite[Chapter IX]{PalWie-B34} and \cite{LotMR97},  which results in PDEs with finite dimensional random inputs.  Hence, our methods are  Wong-Zakai type approximations \cite{WongZak65b,WongZak65},   where  the Brownian motion is approximated by a smooth stochastic process of bounded variation, e.g.,
the spectral approximation used here and piecewise linear approximation of the Brownian motion \cite{WongZak65b,WongZak65}.    We note that piecewise linear approximation can be
used  instead    but this is beyond the scope of the paper.

The resulting PDEs  can be solved numerically using a variety of space-time discretization methods
and any sampling methods or functional expansion methods in random space.
In random space,   we will employ
functional expansion methods,
WCE \cite{BudKal96,LotMR97}
and SCM \cite{ZhangTRK13}, instead of the Monte Carlo method.
These functional expansion methods have no
statistical errors  as no random number generators are used;   they have only errors
from truncations of  Wiener processes and   functional expansions
and allow efficient short-time integration of  SPDEs
 \cite{BudKal96,BudKal97,HouLRZ06,LotMR97,LotRoz06,ZhangRTK12,ZhangTRK13}.

In principle, we can employ any functional expansion, however,
different expansions are preferred for different stochastic products because of computational efficiency.
In practice, WCE is associated  with the Ito-Wick product, see \eqref{eq:sadv-diff-wick}, as the product is
defined with Wiener chaos modes yielding a weakly coupled system (lower-triangular system) of PDEs for linear
equations. On the other hand, SCM is  associated with the Stratonovich product, see \eqref{eq:sadv-diff-stra-BM},
yielding a decoupled
system of PDEs. These different formulations lead to different numerical performance as we demonstrate in Section
\ref{sec:wce-scm-compare}; in particular, WCE can be of second-order convergence in time while SCM is only
of first-order in time
in the second-order moments  for commutative noises.    Further, when the noises serve  as the advection coefficients, SCM can be more
accurate than WCE when both methods are of first order convergence as
the SCM (Stratonovich formulation) can lead to   smaller diffusion coefficients than those for WCE (Ito formulation).

However, a fundamental limitation of these expansion methods is the exponential growth of error with time
and the increasing complexity as the number of random variables is increasing,
generated by the discretization of the Brownian motion.
To deal with this complexity,  a recursive  WCE method was proposed in \cite{LotMR97}
for the Zakai equation
of nonlinear filtering with uncorrelated observations.
More recently, a recursive multi-stage approach was developed to  efficiently  solve  linear stochastic
advection-diffusion-reaction equations using
WCE \cite{ZhangRTK12} or SCM \cite{ZhangTRK13}.

 To deal with the complexity in random space,  some preprocessing  procedures have been
proposed, see e.g. \cite{ChengHZ13,SapLer09}.
In these procedures, we are seeking the solution in the form
$u(t,x;\omega)=  \mean{u(t,x;\cdot)} +\sum_{i=1}^\infty Y_i(t,\omega) u_i(t,x)$.
Then by  imposing the   spatial orthogonality of   $u_i(x,t)$ and  $\partial_t u_j(t,x)$ ($i,j=1,2,\ldots$),
we can obtain an equivalent systems of SPDEs:    a PDE for   $ \mean{u(t,x;\cdot)}$,
a system of equations for $Y_i(t,\omega)$   and  a system of equations  for  $u_i(t,x)$.
In many applications,
this procedure is efficient, even with
few  terms of $Y_i(t,\omega)$ and $u_i(t,x)$ as it may take advantage  of  some intrinsic sparsity structures  of the underlying problems.  However,
this procedure requires  some numerical methods to obtain $Y_i(t,\omega)$, such as WCE (e.g. in \cite{ChengHZ13}) and Monte Carlo methods (e.g. in \cite{SapLer09}).
When WCE or SCM  is used,  the complexity in random space  is still high and thus  the procedure is not efficient for problems driven by   Brownian motion, where many modes $Y_i(t,\omega)$ and $u_i(t,x)$
are required.
Though   these procedures can be applied here,   we  limit ourselves
to the   issue of using the deterministic integration methods, without using these procedures.

Some numerical results of WCE for SPDEs have been presented  in \cite{ZhangRTK12} for linear
advection-diffusion-reaction equations
and in \cite{HouLRZ06} for nonlinear SPDEs including the stochastic Burgers equation
and the Navier-Stokes equations.   Numerical results for SCM have also been provided
 in \cite{ZhangTRK13} for  linear stochastic advection-diffusion-reaction equations
and the stochastic Burgers equation.    Numerical results  have demonstrated that WCE   \cite{ZhangRTK12} and SCM  \cite{ZhangTRK13}  in conjunction with the recursive multi-stage
approach are efficient for long-time integration of
linear advection-diffusion-reaction equations.

The main aim of the current paper is the derivation of theoretical error estimates for both
WCE and SCM methods and subsequent comparison of the numerical
performance of the two methods for commutative
and non-commutative noises.
In addition, we
will develop a recursive multi-stage SCM, different than in \cite{ZhangTRK13}, using a spectral truncation of Brownian motion.
Specifically, in this paper we will derive the error estimate of  WCE for  linear advection-diffusion-reaction
equations
with white noise in the advection velocity  and that of   SCM   with white noise in the reaction rate.
We note that the convergence rate of  WCE is
known only for linear  advection-diffusion-reaction equations with white noise in the reaction rate
although the convergence of  WCE for
linear advection-diffusion-reaction equations has been studied for some time
\cite{LotRoz04,LotMR97,LotRoz06,LotRoz06a}.

The paper is organized as follows. After the Introduction section, in Section \ref{sec:review-scm-wce},  we review the WCE    and SCM for linear parabolic SPDEs and develop a new recursive
SCM using a spectral truncation of Brownian motion, following the same recursive procedure
as WCE in \cite{LotMR97,ZhangRTK12,ZhangTRK13}. In Section \ref{eq:error-estimate-wce-scm} we present the error estimates for both methods
for linear advection-diffusion-reaction equations, with the proofs presented in Section \ref{eq:error-estimate-wce-scm-proof}.
In Section 4, we present numerical results of
WCE and SCM for linear SPDEs with both commutative  and non-commutative noises
and verify the error estimates of WCE and SCM for commutative noises.

 \section{Review of Wiener chaos and stochastic collocation} \label{sec:review-scm-wce}
 In this section, we briefly review WCE and  SCM  for the following linear SPDE in the Ito  form:
\begin{eqnarray}\label{eq:sadv-diff}%
d{u(t,x)}&=&  \mathcal{L}u(t,x)  \,dt+\sum_{k=1}^{q } \mathcal{M}%
_{k}u(t,x)  \,  {d}w_{k}(t),\;(t,x)\in(0,T]\times\cD,\notag \\
{u(0,x)}&=&u_{0}(x),\ \ x\in\cD,
\end{eqnarray}
where   $(w(t),\mathcal{F}_{t})=(\set{ w_{k}(t), 1\leq k\leq q}, \mathcal{F}_{t})$ is a system of one-dimensional independent standard Wiener
processes   defined on a complete  probability space $(\Omega, \cF, P)$   and
\vskip -12pt
\begin{eqnarray} \label{eq:sadv-diff-coefficients}%
\mathcal{L}u(t,x)&=&\sum_{i,j=1}^{d}a_{i,j}\left(  x\right)  D_{i}D_{j}u(t,x)+\sum_{i=1}^{d}b_{i}(x)D_{i}u(t,x)+c\left(  x\right)  u(t,x),\notag\\
\mathcal{M}_{k}u(t,x)&=&\sum_{i=1}^{d}\sigma_{i,k}(x)D_{i}u(t,x)+\nu_{k}\left(
x\right)  u(t,x),
\end{eqnarray}\vskip -8pt
\noindent and $D_{i}$  is the spatial derivative in $x_i$-direction.   We  assume that the domain $\cD$  in $\Real^{d}$  is   such that the periodic boundary conditions can be imposed  or that $\mathcal{D}=\mathbb{R}^{d}.$   
 In the former case, we will consider periodic boundary conditions and in the latter the Cauchy problem.

 We   assume    that  there exist a constant $\delta_{\cL}>0$  and a real number  $C_{\cL}$ such that for any $v\in H^1(\cD)$,
 \begin{equation}\label{eq:determ-strong-parabolic-cond}
\langle\cL v, v \rangle+\frac{1}{2}\sum_{k=1}^q \norm{\cM_k v}^2 + \delta_{\cL}\norm{v}^2_{H^1}\leq C_{\cL}\norm{v}^2,
\end{equation}
 where  $\langle\cdot, \cdot \rangle$ is the duality between the Sobolev spaces $H^{-1}(\cD)$ and $H^1(\cD)$ associated with
the inner-product over $L^2(\cD)$ and  $\norm{\cdot}$ is the    $L^2(\cD)$-norm.
Specifically, we require that    the coefficients of operators
$\mathcal{L}$ and $\mathcal{M}$ in \eqref{eq:sadv-diff-coefficients} are uniformly bounded   and that
\[\sum_{i,j=1}^d \left( 2a_{i,j}(x) -\sum_{k=1}^q\sigma_{i,k}(x)\sigma_{k,j}(x)\right)y_iy_j\geq  2\delta_{\cL}\abs{y}^2,\quad x,y\in \cD,\]
in addition to the Lipschitz continuity of $a_{i,j}(x)$.  If $\mean{\norm{u_0}^2}$ is bounded ($\mean{\cdot}$ is the expectation with respect to $P$),
these assumptions are sufficient for  a unique square-integrable solution of \eqref{eq:sadv-diff}-\eqref{eq:sadv-diff-coefficients}, see e.g.  \cite{LotRoz06a,MikRoz98}.

 The problem  \eqref{eq:sadv-diff}-\eqref{eq:sadv-diff-coefficients} is said to have    commutative noises  if
 \begin{equation} \label{eq:sadv-diff-commutative}
\cM_k \cM_j = \cM_j\cM_k,\quad 1\leq k,j\leq q,
 \end{equation}
 and to have  non-commutative noises otherwise. When $q=1$,  \eqref{eq:sadv-diff-commutative} is satisfied and thus this  is a special case of commutative noises.
 When $\cM_k$ are zeroth-order operators, ($\sigma_{i,k}=0$), \eqref{eq:sadv-diff-commutative} is satisfied and the problem also has commutative noises.
 The definition is  consistent with  that of
  commutative and non-commutative noises for stochastic ordinary differential equations, see e.g. \cite{MilTre-B04}.  In Section \ref{sec:wce-scm-compare}, we  test our algorithms  on  examples  with both  commutative and non-commutative noises.

\begin{rem}
The problem  \eqref{eq:sadv-diff}-\eqref{eq:sadv-diff-coefficients} can be regarded as an approximation of   a problem driven by  a cylindrical Wiener process.
 Consider a cylindrical  Wiener process $W(t,x)=\sum_{k=1}^\infty  \lambda_k  w_k(t)e_k(x)$  where $\sum_{k=1}^\infty\lambda_k^2<\infty$, $\set{w_k(t)}$ are independent Wiener processes,  and
  $\set{e_k(x)}_{k=1}^\infty$ is  a complete orthonormal  basis     (CONS) in $L^2(\cD)$, see e.g. \cite{DaPZab-B92,Roz-B90}. Thus, we can view    \eqref{eq:sadv-diff}-\eqref{eq:sadv-diff-coefficients} as approximations of SPDEs driven by this cylindrical Wiener process.
 \end{rem}

 In both WCE and SCM, we discretize the Brownian motion using the following spectral representations
   (see e.g. \cite{LotMR97,ZhangRTK12}):
    \begin{equation}\label{eq:bm-spectral-exp}
   \lim_{\mathsf{n}\to\infty} \mean{(w(t)-w^{(\mathsf{n})}(t))^2}=0, \quad  w^{(\mathsf{n})}(t)=\sum_{i=1}^\mathsf{n}\xi_{i}
   \int_{0}^{t}m_{i}(s)\,ds,~t\in[0,T],
    \end{equation}
    where  $\xi_{i}$ are mutually
    independent standard Gaussian random variables and $\set{m_{i}}_{i=1}^\infty$ is a   CONS
 in
    $L^2([0,T])$.  The expansion  \eqref{eq:bm-spectral-exp} is
    an extension of Fourier expansion of Brownian motion that is the Wiener construction \cite[Chapter IX]{PalWie-B34},
    see also \cite{KarShr-B91,Kry-B95}.

 \subsection{Wiener chaos expansion (WCE)}
 The WCE solution to \eqref{eq:sadv-diff} is defined   with the Cameron-Martin basis \cite{CamMar47} in Wiener chaos space, using
 Fourier-Hermite series.
The corresponding coefficients   are obtained by solving the
associated  {\em propagator}, which is a
{\em lower-triangular} linear system
of deterministic parabolic equations  determined  by \eqref{eq:sadv-diff}. Specifically, the solution to \eqref{eq:sadv-diff} can be represented as
\begin{equation}\label{eq:sadv-diff-wce-solution}
u(t,x)=\sum_{\alpha\in\cJ_{q }}\frac{1}{\sqrt{\alpha !}}\varphi_{\alpha}(t,x;  u_0)\xi_\alpha, \quad t\in (0,T],
\end{equation}
where   $\cJ_{q }$ is the set
of multi-indices $\alpha=(\alpha_{k,l})_{k,l\geq1}$ of finite length, i.e.,%
\[\cJ_{q }=\left\{  \alpha=(\alpha_{k,l},\ 1\leq k \leq q ),\ l \geq1,\ \ \alpha_{k,l}%
\in\{0,1,2,\ldots\},\ |\alpha|:=\sum_{k,l}\alpha_{k,l}<\infty\right\}.\]
The random variables $\xi_{\alpha}$ are Cameron-Martin orthonormal basis, defined as
\begin{equation}\label{eq:random-basis-temporal-white-noise}%
\xi_{\alpha}:=\prod_{\alpha} \left(  \frac{H_{\alpha_{k,l}}%
(\xi_{k,l})}{\sqrt{\alpha_{k,l}!}}\right)  ,\ \ \ \alpha\in\mathcal{J},
\end{equation}
where
$\displaystyle
\xi_{k,l} =\int_{0}^{T}m_{l}(s)\,dw_{k}(s),
$
and $H_{n}$ is the $n$-th Hermite polynomial:
\begin{equation}\label{eq:Hermite}%
H_{n}(x)=(-1)^{n}e^{x^{2}/2}\frac{d^{n}}{dx^{n}}e^{-x^{2}/2}.
\end{equation}
 Under our assumptions, the SPDE  \eqref{eq:sadv-diff} can   be written in the following form using the Ito-Wick product, see e.g. \cite[Section 2.5]{HolOUZ-B96} and \cite{LotRoz06},
\begin{eqnarray}\label{eq:sadv-diff-wick}%
d{u(t,x)}&=& \mathcal{L}u(t,x) \,dt+\sum_{k=1}^{q }   \mathcal{M}%
_{k}u(t,x)   \diamond \,d{w}_k(t),\;(t,x)\in(0,T]\times\cD,\notag \\
{u(0,x)}&=&u_{0}(x),\ \ x\in\cD,
\end{eqnarray}
where the Ito-Wick product ``$\diamond$''    is defined for the Cameron-Martin basis \eqref{eq:random-basis-temporal-white-noise}  such that
$\displaystyle\xi_{\alpha}\diamond \xi_\beta=\sqrt{\frac{(\alpha+\beta)!}{\alpha!\beta !}}\xi_{\alpha+\beta}.$
 By \eqref{eq:sadv-diff-wce-solution} and Cameron-Martin  theorem \cite{CamMar47}, we  obtain
  the coefficients $\varphi_{\alpha}(t,x;u_0)=\mean{ \sqrt{\alpha !}u(t,x)\xi_\alpha}$.  Approximating $w_k$ with $w_k^{(\mathsf{n})}$ in
\eqref{eq:bm-spectral-exp},   we substitute
  the  representation \eqref{eq:sadv-diff-wce-solution} into  \eqref{eq:sadv-diff-wick}
  and then we can readily check that   the coefficients $\varphi_{\alpha}(t,x; ) $ from \eqref{eq:sadv-diff-wce-solution}  satisfy the
following propagator, see e.g. \cite{LotRoz04,LotRoz06},
\begin{eqnarray}\label{eq:sadv-diff-propagator-general}
\frac{\partial\varphi_{\alpha}(t,x; u_0)}{\partial t}&=&\mathcal{L}%
\varphi_{\alpha}(t,x;u_0)
 +\sum_{k=1}^q \sum_{l=1}^{\mathsf{n}}\alpha_{k,l}m_{l}(s)  \mathcal{M}_{k}\varphi_{\alpha
^{-}(k,l)}(x;  u_0)    ,\quad s\in(0,T],\nonumber\\
\varphi_{\alpha}(0,x)&=&  u_0(x)\mathbf{1}_{\left\{  \left\vert
\alpha\right\vert =0\right\}  },\nonumber
\end{eqnarray}
where
$\alpha^{-}(k,l)$ is the multi-index with components
\begin{equation}
\left(  \alpha^{-}(k,l)\right)_{i,j}=\left\{
\begin{array}
[c]{ll}%
\max(0,{\alpha}_{i,j}-1), & \text{if }i=k\text{ and }j=l,\\
{\alpha}_{i,j}, & \text{otherwise}.
\end{array}
\right.
\end{equation}

 In practical computations, we have to truncate the propagator   \eqref{eq:sadv-diff-propagator-general} and,  consequently, we are interested in
 the following truncated Wiener chaos solution:
\begin{equation}\label{eq:sadv-diff-solution-trun}%
u_{\mathsf{N},\mathsf{n}}(t,x)=\sum_{\alpha\in\cJ_{\mathsf{N},\mathsf{n},q }}\frac{1}{\sqrt{\alpha!}}%
\varphi_{\alpha}(t,x;  u_0)\xi_{\alpha},
\end{equation}
where the set $\cJ_{\mathsf{N},\mathsf{n},q }$
$=\set{\alpha=(\alpha_{k,l})_{q \times \mathsf{n}}|  \sum_{k=1}^q \sum_{l=1}^\mathsf{n}
\alpha_{k,l}\leq \mathsf{N}}$.
Here $\mathsf{N}$ is the highest Hermite polynomial order and $\mathsf{n}$ is the maximum number of
Gaussian random variables for each Wiener process.
In \eqref{eq:random-basis-temporal-white-noise}, we choose the basis $\left\{  m_{l}(s)\right\}  _{l\geq1}$ as
\begin{equation} \label{eq:basis}
m_{1}(s)=\frac{1}{\sqrt{T}},\quad m_{l}(s)=\sqrt{\frac{2}{T}}\cos(\frac
{\pi(l-1)s}{T}),\ \ \ l\geq2,\quad0\leq s\leq T.%
\end{equation}

\begin{figure}[ptb]
\caption{Illustration of the idea of recursive multi-stage approach for long-time integration in \cite{ZhangRTK12}. }%
\label{fig:illus-pic-rm}
\setlength{\unitlength}{0.6mm} \linethickness{1.0pt}
\par
\begin{center}
\begin{picture}(160,30)
\multiput(0,10)(4,0){10}{\textcolor[rgb]{1.00,0.00,0.00}{\line(1,0){1}}}
\put(0,10){\textcolor[rgb]{1.00,0.00,0.00}{\line(0,1){4}}}
\put(2,10){\textcolor[rgb]{1.00,0.00,0.00}{\line(0,1){2}}}
\put(4,10){\textcolor[rgb]{1.00,0.00,0.00}{\line(0,1){2}}}
\put(6,10){\textcolor[rgb]{1.00,0.00,0.00}{\line(0,1){2}}}
\put(8,10){\textcolor[rgb]{1.00,0.00,0.00}{\line(0,1){2}}}
\put(10,14){\textcolor[rgb]{1.00,0.00,0.00}{$\cdots$}}
\put(18,10){\textcolor[rgb]{1.00,0.00,0.00}{\line(0,1){2}}}
\put(20,10){\textcolor[rgb]{1.00,0.00,0.00}{\line(0,1){2}}}
\put(22,10){\textcolor[rgb]{1.00,0.00,0.00}{\line(0,1){2}}}
\put(24,10){\textcolor[rgb]{1.00,0.00,0.00}{\line(0,1){4}}}
\put(48,10){\line(0,1){3}}
\multiput(24,10)(4,0){36}{\line(1,0){3}}
\put(96,10){\line(0,1){3}}
\put(168,10){\line(0,1){3}}
\put(-2,14){$0$}
\put(24,14){$\textcolor[rgb]{1.00,0.00,0.00}{\Delta}$}
\put(50,1){$ $}
\put(48,14){$\textcolor[rgb]{0.00,0.00,1.00}{2\Delta}$}
\put(70,14){\textcolor[rgb]{0.00,0.00,1.00}{$\cdots\cdots$}}
\put(96,14){\textcolor[rgb]{0.00,0.00,1.00}{$i\Delta$}}
\put(120,14){\textcolor[rgb]{0.00,0.00,1.00}{$\cdots\cdots$}}
\put(165,14){$\textcolor[rgb]{0.00,0.00,1.00}{ T=\mathsf{K}\Delta }$}
\put(2,14){\textcolor[rgb]{1.00,0.00,0.00}{$\delta t$} }
\end{picture}
\end{center}
\end{figure}
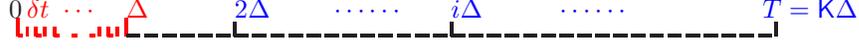

As shown in \cite{BudKal96,LotMR97,ZhangRTK12}, the error  induced by the truncation of Wiener chaos expansion
grows exponentially with time and thus WCE is not efficient for long-time integration.
 To control the error behavior,
\cite{ZhangRTK12} proposes a recursive   WCE (see Algorithm \ref{algo:sadv-diff-s4-mom} below) for  computing the second
moments, $\mean{u^{2}(t,x)}$, of
the solution of the SPDE   \eqref{eq:sadv-diff}.
Specifically, we    discretize  the Brownian motion using the following spectral representation in  a multi-element version, i.e.,
 using $\mathsf{K}$ multi-elements \cite{LotMR97,ZhangRTK12}:
    \begin{equation}\label{eq:multi-elem-spectral-exp}
    w^{(\Delta,\mathsf{n})}(t)=\sum_{k=1}^{\mathsf{K}}\sum_{i=1}^\mathsf{n} \int_{\mathsf{t}_{k-1}\wedge t}^{\mathsf{t}_{k}\wedge t}
    m_{i,k}(s)\,ds\xi_{i,k},~t\in[0, T],
    \end{equation}
    where $0=\mathsf{t}_0<\mathsf{t}_1<\cdots<\mathsf{t}_{\mathsf{K}}=T$, $\mathsf{t}_{k}\wedge t$
    is the minimum of  $\mathsf{t}_k=k\Delta$ and $t$,   $\set{m_{i,k}}_{i=1}^\infty$ is a
    CONS in
    $L^2([\mathsf{t}_k,\mathsf{t}_{k+1}])$, and $\xi_{i,k}$ are mutually
    independent standard Gaussian random variables.
After the truncation of Brownian motion, we can have a similar propagator as \eqref{eq:sadv-diff-propagator-general}. Noticing
the linear property and Markovian properties of the solution to \eqref{eq:sadv-diff},  we take the solution
at $\mathsf{t}_{k-1}$ as initial condition to solve  the solution over $(\mathsf{t}_{k-1},\mathsf{t}_k]$.
Thus, we can recursively compute  the covariance matrix   at   $\mathsf{t}_k$  with the covariance matrix at  the  time instant  $\mathsf{t}_{k-1}$.
We  then have the
 following algorithm for the second moments of the approximate solution; see Figure  \ref{fig:illus-pic-rm} for an illustration and \cite{ZhangRTK12} for the derivation.

 %
\begin{algo}[Recursive multi-stage Wiener chaos expansion, {\cite[Algorithm 2]{ZhangRTK12}}] \label{algo:sadv-diff-s4-mom}
Choose the algorithm's parameters: a CONS $\{e_{m}%
(x)\}_{m\geq1}$ and its truncation $\{e_{m}(x)\}_{m=1}^{\mathsf{M}}$; a time step
$\Delta$; $\mathsf{N}$ and $\mathsf{n}$ which together with the number of noises $q$ determine the size of the
multi-index set $\mathcal{J}_{\mathsf{N},\mathsf{n},q }.$

 Step 1. For each $m=1,\ldots,\mathsf{M},$ solve the propagator
\eqref{eq:sadv-diff-propagator-general} for $\alpha\in\mathcal{J}%
_{\mathsf{N},\mathsf{n},q }$ on the time interval $[0,\Delta]$ with the initial condition
$\phi(x)=e_{m}(x)$ and denote the obtained solution as $\varphi_{\alpha
}(\Delta,x;e_{m}),$  $m=1,\ldots,\mathsf{M},$  and $\alpha\in\cJ_{\mathsf{N},\mathsf{n},q }$
Also, choose  a time step size $\delta t$ to solve numerically the equations in the propagator.

Step 2. Evaluate $q_{\alpha,l,m}=(\varphi_{\alpha}(\Delta,\cdot;e_{l}),e_{m}(\cdot)),$
$l,m=1,\ldots,\mathsf{M}.$   Here  $(\cdot,\cdot)$ is the inner product in $L^2(\cD)$.

Step 3.  Recursively compute the covariance matrices $Q_{lm}({\mathsf{t}_{i}%
};\mathsf{N},\mathsf{n},\mathsf{M})$, $l,m=1,\ldots,\mathsf{M},$ $\mathsf{t}_i=i\Delta$, as follows:
\begin{eqnarray*}
Q_{lm}({0;\mathsf{N},\mathsf{n},\mathsf{M}})  &   =&(u_{0},e_{l})(u_{0},e_{m}),\\
Q_{lm}({\mathsf{t}_{i};\mathsf{N},\mathsf{n},\mathsf{M}})  &  =&\sum_{j,k=1}^{\mathsf{M}}Q_{jk}({\mathsf{t}_{i-1};\mathsf{N},\mathsf{n},\mathsf{M}})\sum_{\alpha\in\mathcal{J}_{\mathsf{N},\mathsf{n},q }%
}\frac{1}{\alpha!}q_{\alpha,j,l}q_{\alpha,k,m},\ i=1,\ldots,\mathsf{K},\nonumber
\end{eqnarray*}
where $u_{0}(x)$ is the initial condition for \eqref{eq:sadv-diff} and obtain $  \mathbb{M}_{\Delta,\mathsf{N},\mathsf{n}}^{\mathsf{M}}(\mathsf{t}_{i},x)$, the second moments
 of the approximate solution to
\eqref{eq:sadv-diff} by the following:
\begin{equation}
 \mathbb{M}_{\Delta,\mathsf{N},\mathsf{n}}^{\mathsf{M}}(\mathsf{t}_{i},x)  =\sum_{l,m=1}^{\mathsf{M}}Q_{lm}({\mathsf{t}_{i};\mathsf{N},\mathsf{n},\mathsf{M}}%
)e_{l}(x)e_{m}(x),\ i=1,\ldots,\mathsf{K}. \label{eq:u2}%
\end{equation}

\end{algo}

\begin{rem}
The complexity of this algorithm is of order $\mathsf{M}^4$ but can be reduced
to the order of $\mathsf{M}^2$ by making full use of the sparsity of the data \cite{ZhangRTK12}.
\end{rem}

 \subsection{Stochastic collocation method (SCM)}
This method leads to a  fully decoupled system
instead of a weakly  coupled system from  the  WCE. First, we rewrite the SPDE
\eqref{eq:sadv-diff} in the Stratonovich form, see e.g. \cite{GyoSti13,Kun82},
    \begin{eqnarray}\label{eq:sadv-diff-stra-BM}%
du(t,x)&=& \tilde{\mathcal{L}}u(t,x) \,dt+\sum_{k=1}^{q } \mathcal{M}_{k}u(t,x)
\circ \,d{w}_k(t),\;(t,x)\in(0,T]\times\cD,\notag \\
u(0,x)&=&u_{0}(x),\ \ x\in\cD,
\end{eqnarray}
where $\tilde{\cL}u=\cL u-\frac{1}{2}\sum_{1\leq k\leq q }\,\mathcal{M}_{k} \mathcal{M}_{k}u$.
Second, we approximate the Brownian motion with its multi-element spectral expansion \eqref{eq:multi-elem-spectral-exp},
and obtain the following partial differential equation with smooth random inputs (see e.g. \cite{GyoSti13}):
    \begin{eqnarray}\label{eq:sadv-diff-trun-BM}%
d{\tilde{u}_{\Delta,\mathsf{n}}(t,x)}&=& \tilde{\mathcal{L}}\tilde{u}_{\Delta,\mathsf{n}}(t,x) \,dt+\sum_{k=1}^q   \mathcal{M}%
_{k}\tilde{u}_{\Delta,\mathsf{n}}(t,x)  {d}w^{(\Delta,\mathsf{n})}_{k}(t),\;(t,x)\in(0,T]\times\mathcal{D},\notag \\
{\tilde{u}(0,x)}&=&u_{0}(x),\ \ x\in\mathcal{D}.
\end{eqnarray}
In \eqref{eq:sadv-diff-trun-BM}, we have   $\mathsf{n}q\mathsf{K}$  standard Gaussian random variables $\xi_{l,k,i}$, $l\leq \mathsf{n},   k \leq q$, $ i\leq  \mathsf{K}$, according to  \eqref{eq:multi-elem-spectral-exp}.
Now we can apply standard numerical techniques of   $\mathsf{n}q\mathsf{K}$-dimensional integration
to numerically obtain $p$-th moments of the solution to  \eqref{eq:sadv-diff-trun-BM}:
 \begin{equation}\label{eq:sadv-diff-truncated-functional}
\mean{\big(
\tilde{u}_{\Delta,\mathsf{n}}(T,x)\big)^p}=\frac{1}{(2\pi)^{\mathsf{n}q\mathsf{K}/2}}\int_{\Real^{\mathsf{n}q\mathsf{K}}}  \big(F(u_0(x),T,x, \mathbf{y}) \big)^p e^{-\frac{\mathbf{y}^\top\mathbf{y}}{2}}
\,d\mathbf{y}, \quad p=1,2,\ldots
\end{equation}
where $\mathbf{y}=(y_{l,k,i})$, $l\leq \mathsf{n},k\leq q$, $i \leq \mathsf{K}$
and the  functional $F$ represents the solution functional for  \eqref{eq:sadv-diff-trun-BM}.
Here we use sparse grid collocation \cite{GenKei96,Smolyak63}
if the dimension $\mathsf{n}q\mathsf{K}$ is moderately large.
As pointed out in \cite{BabNT07,XiuHes05}, we are led to a
fully  {\em decoupled}  system of equations  as in the case of Monte Carlo methods.

In practice, we use the following sparse grid quadrature rule for a $d$-dimensional function $\varphi$, see e.g. \cite{GenKei96,Smolyak63},
\begin{equation}\label{eq:smolyak-tensor-like}
A(\mathsf{L},d)\varphi =\sum_{\mathsf{L}\leq \abs{\mathbf{i} }\leq
\mathsf{L}+d-1}(-1)^{\mathsf{L}+d-1-\left\vert \mathbf{i}\right\vert }\binom{d-1}{\left\vert
\mathbf{i}\right\vert -\mathsf{L}}Q_{i_{1}}\otimes \cdots \otimes  Q_{i_{d}}\varphi,
\end{equation}%
where we have one-dimensional Gauss--Hermite quadrature rules $Q_{n}$ for univariate
functions $\psi (\mathsf{y}),$ $\mathsf{y\in }\mathbb{R}$:
$
Q_{n}\psi (\mathsf{y})=\sum_{\alpha =1}^{n}\psi (\mathsf{y}%
_{n,\alpha })\mathsf{w}_{n,\alpha },$
  $\mathsf{y}_{n,1}<\mathsf{y}_{n,2}<\cdots <\mathsf{y}_{n,n}$ are the
roots of the $n$-th Hermite polynomial \eqref{eq:Hermite}
and $\mathsf{w}_{n,\alpha }$ are the associated weights
$\mathsf{w}_{n,\alpha }=\sfrac{n!}{n^{2}[H_{n-1}(\mathsf{y}_{n,\alpha })]^{2}}.$
The number of sparse grid points, denoted  by $\eta(\mathsf{L},d)$, for this sparse grid quadrature rule is of order
$d^{\mathsf{L}-1}$ when $\mathsf{L}\leq d$, which can be checked readily from the rule \eqref{eq:smolyak-tensor-like}.
For example, we have, for $\mathsf{L}=2,3,4$,
\[\eta(2,d)=2d+1,~\eta(3,d)=2d^2+2d+1,~\eta(4,d)=\frac{4}{3}d^3+2d^2+\frac{14}{3}d+1.\]
 Denote  the set of $\eta(\mathsf{L},d)$ sparse grid points  $ \mathsf{x}_\kappa= ( \mathsf{x}_\kappa^1, \cdots,  \mathsf{x}_\kappa^d) $
 by $\cH_{\mathsf{L}}^{\mathsf{n}q}$, where   $ \mathsf{x}_\kappa^j$ ($1\leq j\leq d$) belongs to the set of points used by the quadrature rule $Q_{i_j}$.
 According to \eqref{eq:smolyak-tensor-like},
we only need to know the function values at the sparse grid
$\cH_{\mathsf{L}}^{\mathsf{n}q}$:
\begin{equation}\label{eq:smolyak-tensor-like-simple}
A(\mathsf{L},d)\varphi =\sum_{\kappa=1}^{\eta(\mathsf{L},\mathsf{n}q) } \varphi(\mathsf{x}_\kappa)\mathsf{W}_\kappa, \quad  \mathsf{x}_\kappa= ( \mathsf{x}_\kappa^1, \cdots,  \mathsf{x}_\kappa^d) \in \cH_{\mathsf{L}}^{\mathsf{n}q},
\end{equation}
where  $\mathsf{W}_\kappa$ are    determined by   \eqref{eq:smolyak-tensor-like} and the choice of the quadrature rules $Q_{i_j}$ and they are called  the sparse grid quadrature weights.

Here again, the direct application of SCM is efficient  only for short-time integration.
To achieve long-time integration, we apply the   recursive multi-stage idea used in Algorithm \ref{algo:sadv-diff-s4-mom}, i.e.,
we  use SCM over
small time interval $(\mathsf{t}_{i-1},\mathsf{t}_i]$
instead of over the whole interval $(0,T]$ and  compute the second-order moments of
the solution recursively in time.  The derivation of such a recursive algorithm  will make use of properties of the problem \eqref{eq:sadv-diff}
and orthogonality of the basis both in
physical space and in random space as will be shown shortly.

We solve \eqref{eq:sadv-diff-trun-BM} with spectral methods in physical space, i.e., using
a truncation of a  CONS in physical space $\set{e_m}_{m=1}^{\mathsf{M}}$
to represent the numerical solution.   The corresponding  approximation  of  $\tilde{u}_{\Delta,\mathsf{n}}(t,x)$
is denoted by $\tilde{u}_{\Delta,\mathsf{n}}^{\mathsf{M}}(t,x)$. Further, let $\upsilon (t,x;s,\upsilon_0)$ be the approximation
$\tilde{u}_{\Delta,\mathsf{n}}^{\mathsf{M}}(t,x)$ of  $\tilde{u}_{\Delta,\mathsf{n}}(t,x)$ with the initial data
$\upsilon_0$ prescribed at $s$:  $\tilde{u}_{\Delta,\mathsf{n}}(s,x)=\upsilon_0(x)$. Note that
\begin{equation}\label{eq:recursive-rep-primitive}
\tilde{u}_{\Delta,\mathsf{n} }^{\mathsf{M}}(\mathsf{t}_i,x)= \upsilon(\mathsf{t}_{i}, x; \mathsf{t}_{i-1},\tilde{u}_{\Delta,\mathsf{n} }^{\mathsf{M}}(\mathsf{t}_{i-1},\cdot)), \quad t_i=i\Delta.
\end{equation}
\noindent Denote $\Phi_m(\mathsf{t}_i; \Delta, \mathsf{n},\mathsf{M})=(\tilde{u}_{\Delta,\mathsf{n}}^{\mathsf{M}}(\mathsf{t}_i,\cdot),e_m)$.   \hskip -2.5pt
Then  the second moments are  computed by
\begin{equation}\label{eq:second-order-moments}
\mean{(\tilde{u}_{\Delta,\mathsf{n}}^{\mathsf{M}}(\mathsf{t}_i,x) )^2}
= \sum_{l,m=1}^{\mathsf{M}}H_{lm}(\mathsf{t}_i;\Delta, \mathsf{n},\mathsf{M})e_l(x)e_m(x),
\end{equation}\vskip -10pt
 \noindent where  $H_{lm}(\mathsf{t}_i; \Delta, \mathsf{n}, \mathsf{M})=\mean{\Phi_l(\mathsf{t}_i; \Delta, \mathsf{n}, \mathsf{M}) \Phi_m(\mathsf{t}_i; \Delta, \mathsf{n},\mathsf{M})}$. \hskip -2pt
Now we show how  the matrix $H_{lm}(\mathsf{t}_i; \Delta, \mathsf{n}, \mathsf{M})$ can be computed recursively.
  By the linearity of    \eqref{eq:sadv-diff-trun-BM}, we have
\[\tilde{u}_{\Delta,\mathsf{n}}^{\mathsf{M}}(\mathsf{t}_i,x)= \sum_{l=1}^{\mathsf{M}}
\Phi_l(\mathsf{t}_{i-1};\Delta,\mathsf{n}, \mathsf{M})
\upsilon(\mathsf{t}_{i},x;\mathsf{t}_{i-1},e_l).\]
 Denote  $h_{l,m,i-1}=(\upsilon(\mathsf{t}_{i},\cdot; \mathsf{t}_{i-1},e_l), e_m)$. Then by the orthonormality of $e_m$,  we  have
\begin{equation*}
\Phi_m(\mathsf{t}_{i}; \Delta,\mathsf{n}, \mathsf{M})= \sum_{l=1}^{\mathsf{M}}
\Phi_l(\mathsf{t}_{i-1}; \Delta, \mathsf{n}, \mathsf{M})h_{l,m,i-1}.
 \end{equation*}
 The matrix $H_{lm}(\mathsf{t}_i; \Delta, \mathsf{n}, \mathsf{M})$ can  be computed recursively as
\begin{equation*}
H_{lm}(\mathsf{t}_i;\Delta,\mathsf{n}, \mathsf{M})=
\sum_{j=1}^{\mathsf{M}}\sum_{k=1}^{\mathsf{M}}H_{jk}(\mathsf{t}_{i-1};\Delta,\mathsf{n}, \mathsf{M})  \mean{h_{j,l,i-1}h_{k,m,i-1} }.
\end{equation*}
We note that the expectation $\mean{h_{j,l,i-1}h_{k,m,i-1} }$ does not depend on $i-1$   because according to  \eqref{eq:sadv-diff-trun-BM} and \eqref{eq:sadv-diff-coefficients},  $\upsilon(\mathsf{t}_{i},x;\mathsf{t}_{i-1},e_l)$ depend  on the length of the time interval $\Delta$  and the random variables $\xi_{l,k,i}$ ($l\leq \mathsf{n}$, $k\leq q$)  but  is independent of time  $\mathsf{t}_{i-1}$.  Denote   $\upsilon(\mathsf{t}_{i},\cdot;\mathsf{t}_{i-1},e_l)$ with $\xi_{l,k,i}$ anchored at the sparse grid point
$\mathsf{x}_\kappa\in  \cH_{\mathsf{L}}^{\mathsf{n}q}$
by $\upsilon_{\kappa}(\Delta,\cdot;e_{l})$.
Let $h_{\kappa,l,m}=(\upsilon_{\kappa}(\Delta,\cdot;e_{l}),e_{m})$.    Then, using the sparse grid quadrature rule   \eqref{eq:smolyak-tensor-like-simple},  we obtain the recursive approximation of  $H_{lm}(\mathsf{t}_i;\Delta,\mathsf{n},\mathsf{M})$:
\begin{equation}\label{eq:recursive-covariance-mat} \resizebox{.90\hsize}{!}
{$\displaystyle
H_{lm}(\mathsf{t}_i;\Delta,\mathsf{n},\mathsf{M})\approx H_{lm}(\mathsf{t}_i;\Delta,\mathsf{L},\mathsf{n},\mathsf{M}) := \sum_{j=1}^{\mathsf{M}} \sum_{k=1}^{\mathsf{M}} H_{jk}(\mathsf{t}_{i-1};\Delta,\mathsf{L},\mathsf{n}, \mathsf{M}) \sum_{\kappa=1}^{\eta(\mathsf{L},\mathsf{n}q)}
h_{\kappa,j,l} h_{\kappa,k,m}\mathsf{W}_{\kappa}.$}
\end{equation}
Substituting  \eqref{eq:recursive-covariance-mat} in  \eqref{eq:second-order-moments},  we obtain  an approximation for the second moments  of $u(t,x)$, denoted by
 $\mathbb{M}_{\Delta, \mathsf{L},\mathsf{n} }^{\mathsf{M}}(\mathsf{t}_i,x)$. When $\mathsf{M}=\infty$ (i.e., when the CONS $\set{e_m}$ is not cut-off), we denote
 this approximation by   $\mathbb{M}_{\Delta, \mathsf{L},\mathsf{n}}(\mathsf{t}_i,x)$.

\begin{rem} For  non-homogeneous equations,  i.e., with  forcing terms, we can have similar
algorithms.  Indeed, the same procedure applies once we can split  the    non-homogeneous equations
into two equations: non-homogeneous equation with zero initial value and  homogeneous equation with initial value.
See \cite{ZhangTRK13} for a derivation of similar algorithms where  only increments of  Brownian motion are used, which is
different from the spectral approximation of Brownian motion used here.
\end{rem}

Now we have the following algorithm for the second moments of the approximate solution.

\begin{algo}[Recursive multi-stage stochastic collocation method]\label{algo:sadv-diff-s4-scm-mom}
Choose  a CONS $\{e_{m}%
(x)\}_{m\geq1}$ and its truncation $\{e_{m}(x)\}_{m=1}^{\mathsf{M}}$; a time step
$\Delta;$ the sparse grid level $\mathsf{L}$ and $\mathsf{n}$, which together with the number of noises $q$ determine
the sparse grid   $\cH_{\mathsf{L}}^{\mathsf{n}q }$  which contains $\eta(\mathsf{L},\mathsf{n}q )$ sparse grid points.

 Step 1. For each $m=1,\ldots,\mathsf{M},$ solve the system of equations
\eqref{eq:sadv-diff-trun-BM} on the sparse grid $\cH_{\mathsf{L}}^{\mathsf{n}q }$ in the time interval $[0,\Delta]$ with the initial condition
$\phi(x)=e_{m}(x)$ and denote the obtained solution as $\upsilon_{\kappa
}(\Delta,x;e_{m}),$  $m=1,\ldots,\mathsf{M},$  and $\kappa=1,\cdots, \eta(\mathsf{L},\mathsf{n}q )$.
Also, choose  a time step size $\delta t$ to solve \eqref{eq:sadv-diff-trun-BM} numerically.

Step 2. Evaluate
$h_{\kappa,l,m}=(\upsilon_{\kappa}(\Delta,\cdot;e_{l}),e_{m}),$
$l,m=1,\ldots,\mathsf{M}.$

Step 3. Recursively compute the covariance matrices $H_{lm}({\mathsf{t}_{i}%
};\mathsf{L},\mathsf{n},\mathsf{M})$, $l,m=1,\ldots,\mathsf{M},$ as follows:
\begin{eqnarray*}
H_{lm}({0;\Delta,\mathsf{L},\mathsf{n},\mathsf{M}})  &   =&(u_{0},e_{l})(u_{0},e_{m}),\\
H_{lm}({\mathsf{t}_{i};\Delta,\mathsf{L},\mathsf{n},\mathsf{M}})  &  =&\sum_{j,k=1}^{\mathsf{M}}H_{jk}({\mathsf{t}_{i-1};\Delta,\mathsf{L},
\mathsf{n}, \mathsf{M}})\sum_{\kappa=1}^{\eta(\mathsf{L},\mathsf{n}q )}%
h_{\kappa,j,l}h_{\kappa,k,m}\mathsf{W}_\kappa,\ i=1,\ldots,\mathsf{K},\nonumber
\end{eqnarray*}
where $u_{0}(x)$ is the initial condition for \eqref{eq:sadv-diff} and obtain the approximate second moments $\mathbb{M}_{\Delta,\mathsf{L},\mathsf{n} }^{\mathsf{M}}(\mathsf{t}_{i},x)$  of the  solution
$u(t,x)$ to
\eqref{eq:sadv-diff} as
\begin{equation}\label{eq:u2-scm}%
 \mathbb{M}_{\Delta,\mathsf{L},\mathsf{n} }^{\mathsf{M} (\mathsf{t}_{i},x)}=\sum_{l,m=1}^{\mathsf{M}}H_{lm}({\mathsf{t}_{i}; \Delta,\mathsf{L},\mathsf{n},\mathsf{M}}%
)e_{l}(x)e_{m}(x),\ i=1,\ldots,\mathsf{K}.
\end{equation}

\end{algo}

\begin{rem}
Similar to Algorithm \ref{algo:sadv-diff-s4-mom}, the cost of this algorithm is $\frac{T}{\Delta}\eta(\mathsf{L}, \mathsf{n}q)\mathsf{M}^4$ and
the storage is  $\eta(\mathsf{L}, \mathsf{n}q)\mathsf{M}^2$. The total cost can be reduced
to the order of $\mathsf{M}^2$ by adopting   reduced order methods
in physical space, see e.g. \cite{SchGit11}. The discussion on computational
efficiency of the recursive WCE methods,  see \cite[Remark 4.1]{ZhangRTK12}, is also valid for Algorithm \ref{algo:sadv-diff-s4-scm-mom}.
\end{rem}


 \section{Error estimates} \label{eq:error-estimate-wce-scm}
Though   WCE  and  SCM  use the same spectral truncation of Brownian motion,
the former is associated with the Ito-Wick product while the latter is related to the Stratonovich product.
Note that
WCE employs orthogonal polynomials   as a basis  and
SCM does not have such orthogonality. This difference allows  WCE to have a better
convergence rate than SCM in the second-order moments, see Corollary  \ref{cor:sadv-diff-cor-estimate-weak}
and Theorem \ref{thm:error-est-wzs-global-local-weak}.

Assume that the operator $\cL$  generates a semi-group $\set{\cT_t}_{t\geq0}$,
which has the following properties:
for $v\in H^k(\cD)$,
\begin{equation}\label{eq:adv-diff-semigroup-estimate}
\norm{\cT_t v}^2_{H^k}\leq  C(k,\cL)e^{2C_{\cL}t}\norm{ v}^2_{H^k},
\end{equation}
where $C(0,\cL)=1$ and
 \begin{equation}\label{eq:adv-diff-semigroup-int-estimate}
\int_{s}^{t}e^{2C_{\cL}(t-\theta)}\norm{\cT_t  v }^2_{H^{k+1}}\,d\theta\leq  \delta_{\cL}^{-1}C(k,\cL)e^{2C_{\cL}(t-s)}\norm{ v}^2_{H^k}.
\end{equation}
 Also, we assume that there exists a constant $\tilde{C}(r,\cM)$ such that
 \begin{equation}\label{eq:adv-diff-estimate-advection-norm}
 \norm{\cM_l  v}_{H^k}^2 \leq
 \tilde{C}(k,\cM)\norm{ v}_{H^{k+1}}^2,\quad \text{for}~ { v}\in H^{k+1},~ l=1,\ldots,q.
 \end{equation}
 and that there exists a constant $\tilde{C}(k,\cL)$ such that
\begin{equation}\label{eq:adv-diff-estimate-diffusion-norm}
\norm{\cL v}_{H^k}^2\leq \tilde{C}(k,\cL) \norm{v}^2_{H^{k+2}}, \text{\, for  }  v\in H^{k+2}.
\end{equation}
The  conditions  \eqref{eq:adv-diff-semigroup-estimate} and \eqref{eq:adv-diff-estimate-advection-norm} are  satisfied  with $k\leq r$
 and  \eqref{eq:adv-diff-estimate-diffusion-norm}  is  satisfied  with $k\leq r-1$
  when  the coefficients from \eqref{eq:sadv-diff-coefficients} belong to the H{\"o}lder space $\cC^{r+1}_b(\cD)$,  which is  equipped with  the following norm
\[\norm{f}_{C^{r}_b}=\max_{0\leq \abs{\beta} \leq \floor{r}} \norm{D_\beta f}_{L^\infty} +
 \sup_{\substack{x,y\in \cD \\ \abs{\beta}=\floor{r}, r>\floor{r}}}\frac{\abs{D_\beta f(x)-D_\beta f(y)}}{\abs{x-y}^{r-\floor{r}}},\]
and $\floor{r}$ is the integer part of the positive   number $r$, cf. \cite[Section 5.1]{Evans-B98}.
 Define  also
 \begin{equation}\label{eq:constant-adv-diff-estimate}
 C_k =\max_{1\leq j\leq k}\set{C(j,\cL)\tilde{C}(j-1,\cM)}.
 \end{equation}

For the  WCE for the SPDE \eqref{eq:sadv-diff}   with  single noise ($q =1$), we have the convergence results stated below.
In the general case, we have not succeeded in proving such theorems but
we numerically check   convergence orders using  examples with commutative noises and non-commutative noises  in   Section  \ref{sec:wce-scm-compare}.
 \begin{thm}\label{thm:error-estimate-propagator-local}
 Let $q=1$ in \eqref{eq:sadv-diff}.
Assume  that $\sigma_{i,1}, a_{i,j}, b_i, c,\nu_1$ in \eqref{eq:sadv-diff-coefficients} belong  to $\cC^{r+1}_b(\cD)$ and $u_0\in H^{r}(\cD)$,
where  $r\geq \mathsf{N}+2$ and $\mathsf{N}$ is the order of Wiener chaos.
 Also assume that \eqref{eq:determ-strong-parabolic-cond}  holds.
Then for $C_1<\delta_{\cL}$, the error of the truncated Wiener  chaos solution  $u_{\mathsf{N},\mathsf{n}}(\mathsf{t}_i,x)$ from   \eqref{eq:sadv-diff-solution-trun} is estimated as
\begin{eqnarray}\label{eq:error-estimate-propagator-local}
&& (\mean{\norm{u_{\mathsf{N},\mathsf{n}}(\mathsf{t}_i,\cdot)-u(\mathsf{t}_i,\cdot)}^2})^{1/2}  \notag\\
 &\leq&
(C_{\floor{r}}\Delta)^{\mathsf{N}/2}e^{C_{\cL}T}\left [\frac{e^{C_{\floor{r}}T}}{(\mathsf{N}+1)!}+ \frac{(C_{\floor{r}}\Delta)^{\floor{r}-\mathsf{N}-1}}{\floor{r}!}
 \frac{\delta_{\cL}}{\delta_{\cL}-C_1}\right]^{1/2}
\norm{u_0}_{H^r}  \notag\\
&&+
\sqrt{2C_{\mathsf{N}+2}C(\mathsf{N}+2,\cL) \tilde{C}(\mathsf{N},\cL)}  e^{C_{\mathsf{N}+2}T+C_\cL T}\frac{\Delta}{\sqrt{\mathsf{n}}\pi}
\norm{u_0}_{H^{\mathsf{N}+2}},
\end{eqnarray}
where  $\mathsf{t}_i=i\Delta$,  the constants $\delta_{\cL}$ and $C_{\cL}$
are from \eqref{eq:determ-strong-parabolic-cond},  $C_{\floor{r}}$ 
 is defined in \eqref{eq:constant-adv-diff-estimate},  $\tilde{C}(\mathsf{N},\cL) $ is from \eqref{eq:adv-diff-estimate-diffusion-norm}, and $C(\mathsf{N}+2,\cL) $ is from \eqref{eq:adv-diff-semigroup-estimate}.
\end{thm}

From Theorem \ref{thm:error-estimate-propagator-local}, we have that the mean-square error of the recursive multi-stage WCE
is $\cO(\Delta^{\mathsf{N}/2})+\cO(\Delta)$.   The same result is proved for $q=1$ and $\sigma_{i,r}=0$ in \cite{LotMR97}, where  the  condition $C_1<\delta_{\cL}$ is not required.   Also, for the case of  $\sigma_{i,r}\neq0$,  the mean-square convergence  without order  but not requiring the condition $C_1<\delta_{\cL}$ was proved in \cite{LotRoz04, LotRoz06}.

\begin{cor}\label{cor:sadv-diff-cor-estimate-weak}
Under the assumptions of Theorem \ref{thm:error-estimate-propagator-local}, we have
\begin{eqnarray}\label{eq:sadv-diff-cor-estimate-weak}
&& \abs{\mean{\norm{u_{\mathsf{N},\mathsf{n}}(\mathsf{t}_i,\cdot)}^2} -\mean{\norm{u(\mathsf{t}_i,\cdot)}^2}}
= \mean{\norm{u_{\mathsf{N},\mathsf{n}}(\mathsf{t}_i,\cdot)-u(\mathsf{t}_i,\cdot)}^2}\notag\\
 &\leq&
(C_{\floor{r}}\Delta)^{\mathsf{N}}e^{2C_{\cL}T}\left [\frac{e^{C_{\floor{r}}T}}{(\mathsf{N}+1)!}+ \frac{(C_{\floor{r}}\Delta)^{\floor{r}-\mathsf{N}-1}}{\floor{r}!}
  \frac{\delta_{\cL}}{\delta_{\cL}-C_1}\right]
\norm{u_0}^2_{H^r} \notag\\
&&+
2C_{\mathsf{N}+2}C (\mathsf{N}+2,\cL) \tilde{C}(\mathsf{N},\cL) e^{2C_{\mathsf{N}+2}T+2C_\cL T}\frac{\Delta^2}{\mathsf{n}\pi^2}
\norm{u_0}^2_{H^{\mathsf{N}+2}}.
\end{eqnarray}
\end{cor}
This corollary   states  that the convergence rate of the error in second-order moments \eqref{eq:sadv-diff-cor-estimate-weak} is twice  that of the mean-square error \eqref{eq:error-estimate-propagator-local} , i.e., $\cO(\Delta^{\mathsf{N}})+\cO(\Delta^2)$.
This corollary can be proved by the orthogonality of  WCE. In fact, it holds that
\begin{equation}\label{eq:wce-orth-second-moment}
\mean{u^2(\mathsf{t}_i,x)}-\mean{u^2_{\mathsf{N},\mathsf{n}}(\mathsf{t}_i,x)}=\mean{
(u(\mathsf{t}_i,x)-u_{\mathsf{N},\mathsf{n}}(\mathsf{t}_i,x))^2},
\end{equation}
as the different terms in the Cameron-Martin basis  are mutually orthogonal \cite{CamMar47}. Then integration
over the physical domain and by the Fubini Theorem, we reach the conclusion in Theorem \ref{thm:error-estimate-propagator-local}.

For SCM for the SPDE \eqref{eq:sadv-diff},  we have the following  estimates: the first one is weak convergence
of the Wong-Zakai type approximation  $\tilde{u}_{\Delta,\mathsf{n}}(t,x)$ from \eqref{eq:sadv-diff-trun-BM} to $u(t,x)$ from \eqref{eq:sadv-diff}, see Theorem
\ref{thm:error-est-wzs-global-local-weak}; 
 the second one is the convergence of SCM, i.e., the convergence of $\mathbb{M}_{\Delta,\mathsf{L},\mathsf{n}}(\mathsf{t}_{i},x)$ to      $\mean{\tilde{u}^2_{\Delta,\mathsf{n}}(\mathsf{t}_i,x)}$, see Theorem
\ref{thm:error-est-wzs-global-local-weak-scm}. Here we
 prove the convergence rate when $\sigma_{i,r}=0$ which belongs to the case of commutative noises \eqref{eq:sadv-diff-commutative}.  Our proof for  Theorem \ref{thm:error-est-wzs-global-local-weak} is based on the  mean-square of convergence of  the Wong-Zakai type approximation \eqref{eq:sadv-diff-trun-BM} to \eqref{eq:sadv-diff}.
 When $\sigma_{i,r}\neq 0$, we  have not succeeded in proving this mean-square convergence and, as far as we know,    only a rate of almost sure convergence of  the Wong-Zakai type approximations to  \eqref{eq:sadv-diff}
 has been proved so far  \cite{GyoSti13}.
\begin{thm}
\label{thm:error-est-wzs-global-local-weak}
Assume  that  $\sigma_{i,r}=0$ and that the initial condition $u_0$ and  the coefficients in \eqref{eq:sadv-diff-coefficients} are in $\cC_b^{2}(\cD)$. Let $u(t,x)$ be the solution to \eqref{eq:sadv-diff} and $\tilde{u}_{\Delta,\mathsf{n}}(t,x)$ be  the solution to
\eqref{eq:sadv-diff-trun-BM}.  Then for any $\varepsilon>0$, there exists a constant  $C>0$ such that  the one-step error  is estimated by
\begin{equation}\label{eq:error-est-wzs-global-local-weak-one-step}
  \abs{\mean{u^2(\Delta,x)}- \mean{\tilde{u}^2_{\Delta,\mathsf{n}}(\Delta,x) }}\leq C\exp(C\Delta) (\Delta^6+\Delta^2)\mathsf{n}^{-1+\varepsilon},
\end{equation}
and the global error   is estimated by
\begin{equation}\label{eq:error-est-wzs-global-local-weak-multi-step}
  \abs{\mean{u^2(\mathsf{t}_i,x)}- \mean{\tilde{u}^2_{\Delta,\mathsf{n}} (\mathsf{t}_i,x)} }\leq C\exp(CT)  \Delta\mathsf{n}^{-1+\varepsilon},\quad 1\leq i\leq\mathsf{K}.
\end{equation}
\end{thm}
 The following theorem is on the convergence of   the second moments  by SCM to those of the solution to \eqref{eq:sadv-diff-trun-BM}. \begin{thm}
\label{thm:error-est-wzs-global-local-weak-scm}
Let   $\tilde{u}_{\Delta,\mathsf{n}}(t,x)$ be the solution to
\eqref{eq:sadv-diff-trun-BM} and  $\mathbb{M}_{\Delta,\mathsf{L},\mathsf{n}}(\mathsf{t}_i,x)$be the limit of  $\mathbb{M}^{\mathsf{M}}_{\Delta,\mathsf{L},\mathsf{n}}(\mathsf{t}_i,x)$ from \eqref{eq:u2-scm} when $\mathsf{M}\to \infty$.  Under the assumptions of Theorem \ref{thm:error-est-wzs-global-local-weak},   for any $\varepsilon>0$,  the one-step error   is estimated by
\begin{eqnarray*}
  \abs{\mathbb{M}_{\Delta,\mathsf{L},\mathsf{n}}(\Delta,x)  - \mean{\tilde{u}^2_{\Delta,\mathsf{n}}(\Delta,x) }}\leq
  C\exp(C\Delta) (\Delta^{3 \mathsf{L}}+\Delta^{2\mathsf{L}}) (1+(3c/2)^{\mathsf{L}\wedge \mathsf{n}})\beta ^{-(\mathsf{L}\wedge \mathsf{n})/2}  \varepsilon^{-\mathsf{L}}{\mathsf{L}}^{-1} \mathsf{n}^{\mathsf{L}\varepsilon},
\end{eqnarray*}
and the global error   is estimated by, for $1\leq i\leq\mathsf{K}$,
\begin{eqnarray*}
  \abs{\mathbb{M}_{\Delta,\mathsf{L},\mathsf{n}}(\mathsf{t}_i,x)- \mean{\tilde{u}^2_{\Delta,\mathsf{n}} (\mathsf{t}_i,x)} }\leq C\exp(CT)  \Delta^{2\mathsf{L}-1} (1+(3c/2)^{\mathsf{L}\wedge \mathsf{n}})\beta ^{-(\mathsf{L}\wedge \mathsf{n})/2}  \varepsilon^{-\mathsf{L}}{\mathsf{L}}^{-1} \mathsf{n}^{\mathsf{L}\varepsilon}.
\end{eqnarray*}
Here  the positive constants  $C,~c,~\beta<1$ are  independent of $\Delta$, $\mathsf{L}$ and $\mathsf{n}$. The expression $\mathsf{L}\wedge \mathsf{n}$ means the minimum of $\mathsf{L}$ and $\mathsf{n}$.
\end{thm}

According to Theorems  \ref{thm:error-est-wzs-global-local-weak} and \ref{thm:error-est-wzs-global-local-weak-scm},
the  error of the SCM is $\cO(\Delta^{2\mathsf{L}-1})+\cO(\Delta)$   in
the second-order moments.
Compared to Corollary  \ref{cor:sadv-diff-cor-estimate-weak},   the SCM is  of one -order lower than   WCE  when $\mathsf{N}=2$ as  the error of
WCE is $\cO(\Delta^{\mathsf{N}})+\cO(\Delta^2).$

 \section{Numerical results}\label{sec:wce-scm-compare}
 In this section, we compare Algorithms \ref{algo:sadv-diff-s4-mom} and \ref{algo:sadv-diff-s4-scm-mom}
 for linear stochastic advection-diffusion-reaction equations with commutative and
 non-commutative noises. We will test the computational performance of these two methods
 in terms of accuracy and computational  cost.
 All the tests  were run using Matlab R2012b, on a Macintosh desktop
computer with Intel Xeon CPU E5462 (quad-core, 2.80 GHz). Every effort was made to program
and execute the different algorithms as much as possible in an identical way.

 We note that we do not have  exact solutions for all examples and hence evaluate the errors of the second-order moments using
 reference solutions, denoted by $\mean{u_{\rm ref}^{2}(T,x)}$, which are obtained by  either Algorithm \ref{algo:sadv-diff-s4-mom} or Algorithm \ref{algo:sadv-diff-s4-scm-mom} with fine resolution.
We do not use solutions obtained from Monte Carlo methods as reference solutions since  Monte Carlo methods are
  of low accuracy 
and are less accurate than the recursive multistage WCE, see \cite{ZhangRTK12} for comparison between WCE and Monte Carlo methods and also below.

The following error measures are  used in  the numerical examples below:
\begin{equation}\label{eq:error-measure-l2-v1}
 \varrho_{2}^{2}(T)   =\abs{\norm{\mean{u_{\rm ref}^{2}(T,\cdot)}}_{l^2}-\norm{\mathbb{M}_{\Delta }^{\mathsf{M}}(T,\cdot)  }_{l^2}},\quad
\varrho_{2}^{r,2}(T) = \frac{\varrho_{2}^{2}(T)}{\norm{\mean{u_{\rm ref}^{2}(T,\cdot)}}_{l^2}},
\end{equation}
\begin{equation} \label{eq:error-measure-linf-v1}
 \varrho_{2}^{\infty}(T)   =\abs{\norm{\mean{u_{\rm ref}^{2}(T,\cdot)}}_{l^\infty}-\norm{\mathbb{M}_{\Delta }^{\mathsf{M}}(T,\cdot))   }_{l^\infty}},\quad
\varrho_{2}^{r,\infty}(T) = \frac{\varrho_{2}^{\infty}(T)}{\norm{\mean{u_{\rm ref}^{2}(T,\cdot)}}_{l^\infty}},
\end{equation}
where $\mathbb{M}^{\mathsf{M}}_{\Delta}(T,x)$ is either  $\mathbb{M}_{\Delta,\mathsf{N},\mathsf{n}}^{\mathsf{M}}(T,x)$ from Algorithm   \ref{algo:sadv-diff-s4-mom}
or $\mathbb{M}_{\Delta,\mathsf{L},\mathsf{n}}^{\mathsf{M}}(T,x)$ from Algorithm \ref{algo:sadv-diff-s4-scm-mom},
$\norm{v}_{l^2}=\displaystyle\left(\frac{2\pi}{\mathsf{M}}\sum_{m=1}^{\mathsf{M}}
  v^2(x_{m})\right)^\frac{1}{2}$, $\norm{v}_{l^\infty}=\displaystyle\max_{1\leq m \leq \mathsf{M}} \abs{v(x_{m})}$, $x_{m}$
  are the Fourier collocation points.

The computational complexity for Algorithm \ref{algo:sadv-diff-s4-mom} is
$\binom{\mathsf{N}+\mathsf{n}q }{\mathsf{N}}\frac{T}{\Delta}\mathsf{M}^4$ (see \cite{ZhangRTK12}) and
that for Algorithm \ref{algo:sadv-diff-s4-scm-mom} is $\eta(\mathsf{L},\mathsf{n}q )\frac{T}{\Delta}\mathsf{M}^4$.
The ratio of the computational cost of SCM  over that of WCE is
$ \sfrac{\eta(\mathsf{L},\mathsf{n}q )}{\binom{\mathsf{N}+\mathsf{n}q }{\mathsf{N}}}$.
 For example, when $\mathsf{N}=1$ and $\mathsf{L=2}$,  the ratio is   $\sfrac{(1+2\mathsf{n}q )}{(1+\mathsf{n}q )}$,
 which will be used in the three numerical examples. The complexity  increases  exponentially with $\mathsf{n}q $ and $\mathsf{L}$, see e.g. \cite{GenKei96}, or $\mathsf{N}$
 but  increases linearly with $\frac{T}{\Delta}$. Hence, we only consider low values of $\mathsf{L}$ and $\mathsf{N}$.

 \begin{exm}[Single noise]\label{exm:sadv-diff-single}
We consider a single noise in  the Ito SPDE \eqref{eq:sadv-diff}   over the domain $(0,T] \times (0,2\pi)$:
\begin{eqnarray}\label{eq:perturbed-sadv-diff}
du &=& [(\epsilon+\frac{1}{2}\sigma^2 )\partial_x^2 u  +  \beta \sin(x)\partial_x u] \,dt+\sigma \partial_x u \,d{w}(t),
\end{eqnarray}
or  equivalently in  the Stratonovich form
\begin{eqnarray}\label{eq:perturbed-sadv-diff-stra}
du &=& [\epsilon\partial_x^2 u  +  \beta \sin(x)\partial_x u] \,dt+\sigma \partial_x u \circ\,d{w}(t),
\end{eqnarray}
with  periodic boundary conditions  and  non-random initial condition $u(0,x)=\cos(x)$,
where $w(t)$ is a standard scalar Wiener process, $\epsilon>0$, $\beta$, $\sigma$ are constants.
\end{exm}

In this example, we compare Algorithms \ref{algo:sadv-diff-s4-mom} and \ref{algo:sadv-diff-s4-scm-mom} for \eqref{eq:perturbed-sadv-diff}
with the parameters  $\beta=0.1$, $\sigma=0.5$ and $\epsilon=0.02$.  We will show that
the recursive multi-stage WCE is at most of order $\Delta^2$ in the second-order moments and the recursive multi-stage SCM is of order $\Delta$.

In Step 1, Algorithm \ref{algo:sadv-diff-s4-mom}, we employ  the   Crank-Nicolson  scheme in time and Fourier collocation
in physical space.
We obtain  the reference solution by  Algorithm \ref{algo:sadv-diff-s4-mom} with the same solver but finer resolution  as a reference solution\footnote{For single noise,
it is proved   in  Theorem \ref{thm:error-estimate-propagator-local} that the  recursive multi-stage WCE is  of second-order convergence in second-order moments. The  second-order convergence is   numerically verified in \cite{ZhangRTK12}.
For this specific example, a Monte Carlo method with $10^6$  sampling paths (which costs $27.6$ hours)  gives  $\norm{\mean{u^2_{\text{MC}}}} = 1.06517\pm 6.1\times10^{-4}$
and  $\norm{\mean{u^2_{\rm MC}}}_\infty = 0.51746\pm6.1\times 10^{-4}$, where the numbers  after `$\pm$' are the statistical errors with the $95\%$ confidence interval. We use  Fourier collocation in space with $M=20$ and Crank-Nicolson in time with  $\delta t=10^{-3}$   for the Monte Carlo method.}
since we have no  exact solution to \eqref{eq:perturbed-sadv-diff}.
The reference solution is obtained by
$\mathsf{M}=30$, $\Delta=10^{-4}$, $\mathsf{N}=4$, $\mathsf{n}=4$, $\delta t=10^{-5}$.
It gives  the second-order moments in $l^2$-norm $\norm{\mean{u^2_{\rm ref}}}_{l^2} \dot{=} 1.065194550063$  and in the $l^\infty$-norm $\norm{\mean{u^2_{\rm ref}}}_{l^\infty}\dot{
=}0.5174746141105$.

From Table   \ref{tbl:multi-stage-wce-advdiff-onenoise-delta-second-order}, we
observe that the recursive  WCE is $\cO(\Delta^{\mathsf{N}})+\cO(\Delta^2)$  for the second-order moments.
When $\mathsf{N}=2$, the method  is  of second-order convergence in $\Delta$ and of first-order convergence when $\mathsf{N}=1$.
When $\mathsf{N}=3$, the method is still second-order in $\Delta$ (not  presented here).
This verifies the estimate  in Corollary \ref{cor:sadv-diff-cor-estimate-weak}.

 \begin{table}[htbp] \centering
  \tabcaption{Algorithm \ref{algo:sadv-diff-s4-mom}: recursive multi-stage Wiener chaos method for
  \eqref{eq:perturbed-sadv-diff} at $T=5$: $\sigma=0.5$, $\beta=0.1$, $\epsilon=0.02$,   and $\mathsf{M}=20$, $\mathsf{n}=1$. }
  \label{tbl:multi-stage-wce-advdiff-onenoise-delta-second-order}
  \begin{tabular}[c]{ccccc|ccc}\hline
     $\Delta$  &    $\delta t$ & $\mathsf{N}$     &   $\varrho_2^{r,2}(T)$    & order     &  $\varrho_2^{r,\infty}(T)$ & order  & CPU time (sec.)   \\    \hline  
 1.0e-1  & 1.0e-2    & 1   & 1.5249e-2  &   --               & 8.8177e-3  &                  & 3.57   \\   \hline
 1.0e-2  & 1.0e-3    & 1   & 1.5865e-3  &  $\Delta^{0.98}$   & 8.9310e-4  & $\Delta^{0.99}$  & 33.22   \\    \hline
 1.0e-3  & 1.0e-4    & 1   & 1.5934e-4  &  $\Delta^{1.00}$   & 8.9429e-5  & $\Delta^{1.00}$  & 348.03   \\    \hline  \hline
 1.0e-1  & 1.0e-2    & 2   & 1.9070e-4  &  --                & 4.1855e-5  & --               & 5.14    \\   \hline
 1.0e-2  & 1.0e-3    & 2   & 2.0088e-6  &  $\Delta^{1.98}$   & 4.2889e-7  & $\Delta^{1.99}$  & 51.75   \\    \hline
 1.0e-3  & 1.0e-4    & 2   & 2.0386e-8  &  $\Delta^{1.99}$   & 4.8703e-9  & $\Delta^{1.94}$  & 490.04   \\    \hline  \hline
 \end{tabular}%
  \end{table}

In Step 1, Algorithm \ref{algo:sadv-diff-s4-scm-mom}, we  use the   Crank-Nicolson scheme  in time and Fourier collocation
method in physical space.  The errors are also measured  as in
\eqref{eq:error-measure-l2-v1} and \eqref{eq:error-measure-linf-v1}.
The reference solution is obtained  by Algorithm  \ref{algo:sadv-diff-s4-mom} as in the case of WCE.
We observe in Table  \ref{tbl:multi-stage-scm-advdiff-onenoise-delta-1st-order}  that the convergence order for second-order moments is one in $\Delta$ even when the sparse grid level
$\mathsf{L}$ is $2$, $3$ and $4$ (the last is not presented here).  The errors for $\mathsf{L}=3$ are more than half  in magnitude smaller than those for $\mathsf{L}=2$ while
the time cost for $\mathsf{L}=3$ is about $1.5$ times of that for $\mathsf{L}=2$.

\begin{table}[htbp] \centering
      \tabcaption{Algorithm \ref{algo:sadv-diff-s4-scm-mom}: recursive multi-stage  stochastic   collocation  method for
  \eqref{eq:perturbed-sadv-diff} at $T=5$: $\sigma=0.5$, $\beta=0.1$, $\epsilon=0.02$,   and $\mathsf{M}=20$, $\mathsf{n}=1$. }
  \label{tbl:multi-stage-scm-advdiff-onenoise-delta-1st-order}
  \begin{tabular}[c]{ccccc|ccc}\hline
     $\Delta$  &    $\delta t$ & $\mathsf{L}$     &   $\varrho_2^{r,2}(T)$    & order     &  $\varrho_2^{r,\infty}(T)$ & order  & CPU time (sec.)   \\    \hline     \hline
   1e-01  & 1e-02   & 2    &3.4808e-04   & --                              & 3.0383e-03  &  --                             &3.71  \\    \hline
   1e-02  & 1e-03   & 2    &3.4839e-05   &$\Delta^{1.00}$   & 3.0130e-04  &$\Delta^{1.00}$  &33.88  \\    \hline
   1e-03  & 1e-04   & 2    &3.4844e-06  &$\Delta^{1.00}$   & 3.0106e-05  & $\Delta^{1.00}$  &325.06 \\    \hline  \hline
   1e-01  & 1e-02   & 3    &1.6815e-04  & --                              & 3.4829e-04  &  --                             & 5.16  \\    \hline
   1e-02  & 1e-03   & 3    &1.6230e-05  &$\Delta^{1.02}$  & 3.2283e-05  & $\Delta^{1.03}$  & 50.59  \\    \hline
   1e-03  & 1e-04   & 3    &1.6170e-06  &$\Delta^{1.00}$  & 3.2026e-06  &$\Delta^{1.00}$   & 486.08 \\    \hline  \hline
 \end{tabular}%

  \end{table}

In summary, from Tables \ref{tbl:multi-stage-wce-advdiff-onenoise-delta-second-order} and  \ref{tbl:multi-stage-scm-advdiff-onenoise-delta-1st-order},
we observe that the recursive multi-stage WCE is $\cO(\Delta^{\mathsf{N}})+\cO(\Delta^2)$ and the recursive multi-stage SCM is $\cO(\Delta)$, as predicted
by the error estimates in Section \ref{eq:error-estimate-wce-scm}.
While  the SCM and   the  WCE
are of the same order when $\mathsf{N}=1$ and $\mathsf{L}\geq2$,  the former can be more accurate than the latter.
In fact, when $\mathsf{N}=1$ and $\mathsf{L}=2$, the recursive multi-stage SCM error is almost two orders of magnitude smaller than
the recursive multi-stage WCE while the computational cost for both is almost the same, as predicted
($\binom{\mathsf{N}+\mathsf{n}q }{\mathsf{N}}=\eta(\mathsf{L},\mathsf{n}q )=2$).
The recursive multi-stage WCE with $\mathsf{N}=2$ is  of order $\Delta^2$ and its errors are almost two orders of magnitude smaller than those by
the recursive multi-stage SCM (with level 2 or 3)  for
the second-order moments.

In this example,  the recursive multi-stage SCM outperforms  the  recursive multi-stage WCE with $\mathsf{N}=1$. The reason can be as follows.
In SCM, we solve an  advection-dominated  equation rather than a  diffusion-dominated equation in WCE, as
SCM is associated with the Stratonovich product which leads to the removal of the term $\frac{1}{2}\sigma^2\partial_x^2u$ in
the resulting equation, see \eqref{eq:perturbed-sadv-diff-stra}.
The larger $\sigma$ is, the more  dominant the diffusion is. In fact,
results for $\sigma=1$ and $\sigma=0.1$ (not presented here) show that
when $\sigma=1$, the relative error of SCM with $\mathsf{L}=2$ is almost three orders of magnitude smaller than WCE with $\mathsf{N}=1$;
when $\sigma=0.1$, the relative error of SCM with $\mathsf{L}=2$ is only less than one order of magnitude smaller than WCE with $\mathsf{N}=1$.
With the Crank-Nicolson  scheme in time and Fourier collocation in physical space, we cannot achieve better accuracy for WCE with $\mathsf{N}=1$ and
$\Delta$ no less than $0.0005$ when $\mathsf{M}\leq 40$.

%
%
 \begin{exm}[Commutative noises]\label{exm:sadv-diff-com}
 We consider two commutative noises in the Ito SPDE \eqref{eq:sadv-diff}   over the domain $(0,T] \times (0,2\pi)$:
\begin{eqnarray}\label{eq:perturbed-sadv-diff-com}
du &=& [(\epsilon+\frac{1}{2}\sigma_1^2\cos^2(x) )\partial_x^2 u  +  (\beta \sin(x)-\frac{1}{4}\sigma_1^2\sin(2x))\partial_x u ]\,dt \notag\\
   & &+\sigma_1 \cos(x)\partial_x u \,dw_1(t) + \sigma_2  u \,dw_2(t),
\end{eqnarray}
or equivalently in the Stratonovich form
\begin{eqnarray}\label{eq:perturbed-sadv-diff-com-stra}
du &=& [\epsilon \partial_x^2 u  +  \beta \sin(x)\partial_x u ]\,dt+\sigma_1 \cos(x)\partial_x u \circ\,dw_1(t) +
\sigma_2  u \circ\,dw_2(t),
\end{eqnarray}
with  periodic boundary conditions and non-random initial condition $u(0,x) = \cos(x)$,
where $(w_1(t),w_2(t))$ is a standard  two-dimensional Wiener process, $\epsilon>0$, $\beta$, $\sigma_1,~\sigma_2$ are constants.
 The problem has  commutative noises, see  \eqref{eq:sadv-diff-commutative}.
\end{exm}

In this example, we take  $\sigma_1=0.5$, $\sigma_2=0.2$, $\beta=0.1$, $\epsilon=0.02$.
We again observe first-order convergence for SCM and WCE with $\mathsf{N}=1$, and second-order convergence for WCE with $\mathsf{N}=2$
as in the last example with single noise.

We choose the same space-time solver for the recursive multi-stage WCE and SCM as in the last example.   We compute the errors  as  in \eqref{eq:error-measure-l2-v1} and \eqref{eq:error-measure-linf-v1}.
In Table   \ref{tbl:multi-stage-wce-advdiff-commut-delta-1st-order},   the reference second moments  are  $\norm{\mathbb{M}^{\mathsf{M}}_{\Delta=10^{-4}, \mathsf{N},\mathsf{n}}(T,\cdot))  }_{l^2}$
 and  $\norm{\mathbb{M}^{\mathsf{M}}_{\Delta=10^{-4}, \mathsf{N},\mathsf{n}}(T,\cdot) }_{l^\infty}$  obtained by  Algorithm  \ref{algo:sadv-diff-s4-mom} with   $\delta t= 10^{-5}$  and
all the other truncation parameters are   the same as stated in the table.
In Table \ref{tbl:multi-stage-scm-advdiff-commut-delta-half-order},
the reference second moments  are $\norm{\mathbb{M}^{\mathsf{M}}_{\Delta=10^{-4},\mathsf{L},\mathsf{n}}(T)}_{l^2}$ and $\norm{\mathbb{M}^{\mathsf{M}}_{\Delta=10^{-4},\mathsf{L},\mathsf{n}}(T)}_{l^\infty}$  obtained by Algorithm   \ref{algo:sadv-diff-s4-scm-mom}  with    $\delta t= 10^{-5}$  while all the other truncation parameters are the same as in the table.

Here we do not compare the performance of Monte Carlo simulations with our algorithms as the main cost of Monte Carlo methods is to reduce the statistical errors.  For the same parameters described above,
when we  used $ 10^6$ Monte Carlo  sampling paths, we could only reach the statistical error  of  $8.3\times 10^{-4}$, in $3.9$ hours. To obtain an error of    $1\times 10^{-5}$,   seven thousands times more
Monte Carlo  sampling paths should be used, requiring three years of computational time and thus  is not considered here.  In the next example, we have similar situations and hence we will not consider Monte Carlo simulations.     This also demonstrates the computational efficiency of Algorithms \ref{algo:sadv-diff-s4-mom} and \ref{algo:sadv-diff-s4-scm-mom} in comparison with Monte Carlo methods.

 For WCE, we observe in  Table \ref{tbl:multi-stage-wce-advdiff-commut-delta-1st-order} convergence of  order $\Delta^{\mathsf{N}}$ ($\mathsf{N}\leq2$) in
 the second-order moments: first-order convergence  when  $\mathsf{N}=1$,
 and second-order convergence when $\mathsf{N}=2$.  Numerical results  for $\mathsf{N}=3$ (not presented here) show that
 the convergence order is still two even though the accuracy is further improved when $\mathsf{N}$ increases from 2 to 3.
 This is consistent with our estimate $\cO(\Delta^{\mathsf{N}})+\cO(\Delta^2)$ in Corollary \ref{cor:sadv-diff-cor-estimate-weak}.

  We also tested the case $\mathsf{n}=2$ which gives similar results and the same convergence order.
\begin{table}[htbp] \centering
  \tabcaption{Algorithm \ref{algo:sadv-diff-s4-mom}: recursive multi-stage Wiener chaos expansion for  commutative noises
  \eqref{eq:perturbed-sadv-diff-com} at $T=1$: $\sigma_1=0.5$, $\sigma_2=0.2$, $\beta=0.1$, $\epsilon=0.02$,
  and $\mathsf{M}=30$, $\mathsf{n}=1$. }
  \label{tbl:multi-stage-wce-advdiff-commut-delta-1st-order}
 \begin{tabular}[c]{ccccc|ccc}\hline
     $\Delta$  &    $\delta t$ & $\mathsf{N}$     &   $\bar{\varrho}_2^{r,2}(T)$    & order     &  $\bar{\varrho}_2^{r,\infty}(T)$ & order  & CPU time (sec.)  \\    \hline     \hline
  1e-01   & 1e-02     &  1   &   1.6994e-03  & --               &  1.6548e-03                   & --                & 3.19      \\    \hline
  1e-02   & 1e-03     &  1   &   1.7838e-04  & $\Delta^{0.98}$  &  1.7172e-04   & $\Delta^{ 0.98}$  & 32.74      \\    \hline
  1e-03   & 1e-04     &  1   &   1.6323e-05  & $\Delta^{1.04}$  &  1.5694e-05   & $\Delta^{ 1.04}$  & 329.15     \\    \hline  \hline
  1e-01   & 1e-02     &  2   &   4.0658e-05   & --               &           2.9568e-05  &      --           & 6.53       \\    \hline
  1e-02   & 1e-03     &  2   &   4.4805e-07   &$\Delta^{1.96}$   &  3.3106e-07  & $\Delta^{ 1.95 }$ & 65.89      \\    \hline
  1e-03   & 1e-04     &  2   &   4.4682e-09  &$\Delta^{2.00}$   &  3.3484e-09   & $\Delta^{ 2.00 }$ & 657.55     \\    \hline      \hline
 \end{tabular}%
\end{table}

  For SCM,  we observe  first-order convergence in $\Delta$ from Table \ref{tbl:multi-stage-scm-advdiff-commut-delta-half-order}
  when $\mathsf{L}=2,3$.
  We note that further refinement in  truncation parameters in random space, i.e., increasing $\mathsf{L}$ and/or $\mathsf{n}$ do not
  change the convergence order nor improve the accuracy.
  The case $\mathsf{L}=3$ actually
  leads to a bit worse accuracy, compared with the case $\mathsf{L}=2$. We tested the case $\mathsf{L}=4$ which leads to  the same magnitudes of errors  as $\mathsf{L}=3$.  We also  tested  $\mathsf{n}=2$ and observed no improved
  accuracy for $\mathsf{L}=2,3,4$. These numerical results are not presented here.

   \begin{table}
   \tabcaption{Algorithm \ref{algo:sadv-diff-s4-scm-mom}: recursive multi-stage  stochastic   collocation method for  commutative noises   \eqref{eq:perturbed-sadv-diff-com} at $T=1$: $\sigma_1=0.5$, $\sigma_2=0.2$, $\beta=0.1$, $\epsilon=0.02$,   and $\mathsf{M}=30$, $\mathsf{n}=1$. }
  \label{tbl:multi-stage-scm-advdiff-commut-delta-half-order}
  \begin{tabular}[c]{ccccc|ccc}\hline
     $\Delta$  &    $\delta t$ & $\mathsf{L}$     &   $\bar{\varrho}_2^{r,2}(T)$    & order     &  $\bar{\varrho}_2^{r,\infty}(T)$ & order  & CPU time (sec.)  \\    \hline     \hline

 1e-01 & 1e-02  &  2   &  1.3624e-04 & --              &1.2453e-03  &   --               & 5.18    \\    \hline
 1e-02 & 1e-03  &  2   &  1.3064e-05 &$\Delta^{1.02}$  &1.2009e-04  &$\Delta^{1.02}$     & 54.70    \\    \hline
 1e-03 & 1e-04  &  2   &  1.1837e-06 &$\Delta^{1.04}$  &1.0889e-05  &$\Delta^{1.04}$     & 545.20   \\    \hline  \hline
 1e-01 & 1e-02  &  3   &  2.5946e-04 & --                              &2.0482e-04  &  --                & 13.26    \\    \hline
 1e-02 & 1e-03  &  3   &  2.5437e-05  &$\Delta^{1.01}$  &1.7897e-05    &$\Delta^{1.06}$     & 142.23   \\    \hline
 1e-03 & 1e-04  &  3   &  2.3102e-06 &$\Delta^{1.04}$   & 1.6062e-06 &$\Delta^{1.05}$     & 1420.24   \\    \hline  \hline

  \end{tabular}%
  \end{table}

For the two commutative noises, we conclude from this example that the recursive multi-stage WCE is of order $\Delta^{\mathsf{N}}+\Delta^2$  in the second-order moments
and that the  recursive multi-stage SCM is of order $\Delta$ in the second-order moments no matter what sparse grid level is used.
The errors of recursive multi-stage SCM are one order of magnitude smaller than those of recursive multi-stage WCE  with $\mathsf{N}=1$ while  the time cost of SCM
is about 1.6 times of that cost of WCE. For large magnitude of noises
($\sigma_1=\sigma_2=1$, numerical results  are not presented),
we observed that the SCM with $\mathsf{L}=2$ and WCE with $\mathsf{N}=1$ have the same order-of-magnitude accuracy.
In this example,  the use of SCM with $\mathsf{L}=2$ for small magnitude of noises
is competitive with the use of WCE with $\mathsf{N}=1$.
%
\begin{exm}[Non-commutative noises]\label{exm:sadv-diff-noncom}
We consider two non-commutative noises in the Ito SPDE \eqref{eq:sadv-diff}   over the domain $(0,T] \times (0,2\pi)$:
\begin{eqnarray}\label{eq:perturbed-sadv-diff-noncom}
du &=& [(\epsilon+\frac{1}{2}\sigma_1^2 )\partial_x^2 u  +  \beta \sin(x)\partial_x u+\frac{1}{2}\sigma_2^2\cos^2(x)u ]\,dt \notag\\
   &&+\sigma_1 \partial_x u \,dw_1(t) +  \sigma_2 \cos(x) u \,dw_2(t),
\end{eqnarray}
or  equivalently in the Stratonovich form
\begin{eqnarray}\label{eq:perturbed-sadv-diff-noncom-stra}
du &=& [\epsilon \partial_x^2 u  +  \beta \sin(x)\partial_x u ]\,dt+\sigma_1 \partial_x u \circ\,dw_1(t) +
\sigma_2 \cos(x) u \circ\,dw_2(t),
\end{eqnarray}
with  periodic boundary conditions and   non-random initial condition $u(0,x) = \cos(x)$,
where $(w_1(t),w_2(t))$ is a standard Wiener process, $\epsilon>0$, $\beta$, $\sigma_1$, $\sigma_2$ are constants.
 The problem has    non-commutative noises as the coefficients do not satisfy  \eqref{eq:sadv-diff-commutative}.
\end{exm}

We take the same constants $\epsilon>0$, $\beta$, $\sigma_1$, $\sigma_2$ as   in the last example.
We also take the same space-time solver as in the last example. In the current example,
we  observe only first-order convergence for SCM (level $\mathsf{L}=2,3,4$) and WCE ($\mathsf{N}=1,2,3$)
when $\mathsf{n}=1,2$, see Table \ref{tbl:multi-stage-wce-advdiff-noncommut-delta-half-order-n1} for parts of
the numerical results.

The errors are computed  as in the last example.  The  reference solutions  are obtained by
Algorithm  \ref{algo:sadv-diff-s4-mom}  for  the  recursive multi-stage WCE  solutions and
by Algorithm  \ref{algo:sadv-diff-s4-scm-mom} for the recursive multi-stage SCM solutions,  with   $\Delta= 5\times 10^{-4}$ and  $\delta t= 5\times 10^{-5}$
 and all the other truncation parameters  the same as    stated in Tables   \ref{tbl:multi-stage-wce-advdiff-noncommut-delta-half-order-n1}
 and  \ref{tbl:multi-stage-wce-advdiff-noncommut-delta-half-order-n10}.

 \begin{table}[htbp] \centering
 \tabcaption{Algorithm \ref{algo:sadv-diff-s4-mom} (recursive multi-stage Wiener chaos expansion, left) and
Algorithm \ref{algo:sadv-diff-s4-scm-mom} (recursive multi-stage stochastic collocation method, right) for
 \eqref{eq:perturbed-sadv-diff-noncom} at $T=1$: $\sigma_1=0.5$, $\sigma_2=0.2$, $\beta=0.1$, $\epsilon=0.02$,
 and $\mathsf{M}=20$, $\mathsf{n}=1$. The time step size $\delta t$ is $\sfrac{\Delta}{10}$. The reported CPU time
 is in seconds.}
 \label{tbl:multi-stage-wce-advdiff-noncommut-delta-half-order-n1}
\scalebox{1.0}{
   \begin{tabular}[c]{ccccc|cccccc}\hline                  \hline
  $\Delta$  &     $\mathsf{N}$     &   $\bar{\varrho}_2^{r,2}(T)$    & order   &   time (sec.)  & $\mathsf{L}$
  &   $\bar{\varrho}_2^{r,2}(T)$    & order    &  time (sec.)\\ \hline \hline
  1.0e-1  & 1   &  3.7516e-03 &  --                           & 1.04               &2   &  6.4343e-04 & --              &1.65     \\    \hline
  5.0e-2  & 1   &  1.8938e-03                    & $\Delta^{0.99}$ & 2.11             &2    &  3.1738e-04 &$\Delta^{1.02 }$  &3.31     \\    \hline
  2.0e-2  & 1   &  7.5292e-04                  & $\Delta^{1.01}$ & 5.12             &2    &   1.2440e-04 &$\Delta^{1.02  }$  &8.64     \\    \hline
  1.0e-2  & 1   &  3.6796e-04                  & $\Delta^{1.03}$ & 10.19            &2   &\textbf{6.0502e-05}   &$\Delta^{1.04 }$  &\textbf{17.12}     \\    \hline
  5.0e-3  & 1   &  1.7457e-04                  & $\Delta^{1.08}$ & 20.01            &2   & 2.8635e-05    &$\Delta^{1.08 }$  &33.82   \\    \hline
  2.0e-3  & 1   &  \textbf{5.8246e-05} & $\Delta^{1.20}$ & \textbf{50.39}            &2   &  9.5401e-06                      &$\Delta^{1.12}$  &86.44    \\    \hline  \hline
%
  1.0e-1  & 2   &   9.4415e-05    &  --                                                     & 2.16                       &3   &  1.5803e-04   & --                 &4.03     \\    \hline
  5.0e-2  & 2   & \textbf{3.7303e-05} & $\Delta^{1.81}$             & \textbf{4.11}        &3   &  \textbf{7.6548e-05}  &$\Delta^{1.05 }$  &\textbf{8.68}     \\    \hline
  2.0e-2  & 2   &    1.2282e-05 & $\Delta^{1.34}$ 	& 9.97              &3   & 2.9673e-05  &$\Delta^{1.03 }$  &22.08    \\    \hline
  1.0e-2  & 2   &    5.5807e-06 & $\Delta^{1.21}$ 	& 20.03             &3   &  1.4378e-05  &$\Delta^{1.05 }$  &43.85    \\    \hline
  5.0e-3  & 2   &     2.5471e-06 & $\Delta^{1.14}$	 & 40.25             &3   &  6.7925e-06  &$\Delta^{1.08 }$  &88.35    \\    \hline
  2.0e-3  & 2   &     8.2965e-07 & $\Delta^{1.22}$ 	& 101.34            &3   &   2.2605e-06  &$\Delta^{1.20 }$  &223.15   \\    \hline  \hline
 %
 \end{tabular}
 }
 \end{table}

 In this example, our error estimate for recursive multi-stage WCE is not valid any more and the numerical results suggest  that
 the errors behave as
 $ \Delta^{\mathsf{N}}+C\sfrac{\Delta}{\mathsf{n}}$.  For $\mathsf{N}=1$ and $\mathsf{n}=10$ (not presented), the error is almost the same
 as $\mathsf{n}=1$.   While $\mathsf{N}=2$ and $\mathsf{n}=10$, the error first decreases as $\cO(\Delta^2)$ for large time step size
 and then as $\cO(\Delta)$ for small time step size; see Table
 \ref{tbl:multi-stage-wce-advdiff-noncommut-delta-half-order-n10}.
 When $\mathsf{N}=2$ and $\mathsf{n}=10$, the errors with $\Delta=0.005,0.002,0.001$ are ten percent ($\sfrac{1}{\mathsf{n}}$) of those
 with the same parameters but $\mathsf{n}=1$ in Table \ref{tbl:multi-stage-wce-advdiff-noncommut-delta-half-order-n1}.
 Here the constant in front of $\Delta$, $C/\mathsf{n}$, plays  an important role:
 when $\Delta$ is large and this constant is small, then the order of two can be observed; when $\Delta$ is small,
 $C\sfrac{\Delta}{\mathsf{n}}$ is dominant so that only first-order convergence can be observed.

 The recursive multi-stage SCM is of first-order convergence when $\mathsf{L}=2,3,4$ and $\mathsf{n}=1,2,10$ (only parts of
 the results presented).    In contrast to  Example \ref{exm:sadv-diff-com},
 the errors from $\mathsf{L}=3$ are one order of magnitude smaller those from $\mathsf{L}=2$.
 Recalling that the number of sparse grid points is $\eta(2,2)=5$ and $\eta(3,2)=13$, we have   the cost for $\mathsf{L}=3$ is about 2.6 times of that for
 $\mathsf{L}=2$.  However, it is expected that in practice, a low level sparse grid is more efficient
 than a high level one when $\mathsf{n}q $ is large
 as the number of sparse grid points  $\eta(\mathsf{L},\mathsf{n}q )$ is increasing exponentially with $\mathsf{n}q $
 and $\mathsf{L}$.  In other words, $\mathsf{L}=2$ is preferred when SPDEs with many noises (large $q$) are considered.

 As discussed  in the beginning of this section,
 the ratio of time cost for SCM and WCE is $\eta(\mathsf{L},\mathsf{n}q )/\binom{\mathsf{N}+\mathsf{n}q }{\mathsf{N}}$.
 The cost of recursive multi-stage SCM with $\mathsf{L}=2$ is   at most 1.8 times  (1.6 predicted by the ratio above, $q =2$ and $\mathsf{n}=1$) of that  of recursive multi-stage WCE with $\mathsf{N}=1$.
 However, in this example, the accuracy of the recursive multi-stage SCM is one order of magnitude smaller  than that of the recursive multi-stage WCE when
$\mathsf{N}=1$ and $\mathsf{L}=2$.
In Table \ref{tbl:multi-stage-wce-advdiff-noncommut-delta-half-order-n1}, we present in bold the errors between $3.5\times 10^{-5}$ and $8.0\times 10^{-5}$.
 Among the four cases listed in the table, the most efficient, for the given accuracy above,    is  WCE with $\mathsf{N}=2$, which outperforms
 SCM with $\mathsf{L}=3$ and $\mathsf{L}=2$.    WCE  with $\mathsf{N}=1$ is less efficient than the other three cases.
 We also observed that  when $\sigma_1=\sigma_2=1$,  SCM  with $\mathsf{L}=2$ is one order of magnitude smaller than
 WCE with  $\mathsf{N}=1$ (results not presented here).

 For non-commutative noises in this example, we show that the error for WCE is $\Delta^2+ C\sfrac{\Delta}{\mathsf{n}}$ and the error for SCM is $\Delta$.
 The numerical results suggest that SCM with $\mathsf{L}=2$ is competitive with WCE with $\mathsf{N}=1$ for both  small and large magnitude of noises
 if $\mathsf{n}=1$.

\begin{table}[htbp] \centering
  \tabcaption{Algorithm \ref{algo:sadv-diff-s4-mom}: recursive multi-stage Wiener  chaos expansion for
  \eqref{eq:perturbed-sadv-diff-noncom} at $T=1$: $\sigma_1=0.5$, $\sigma_2=0.2$, $\beta=0.1$, $\epsilon=0.02$. The parameters are
  $\mathsf{M}=20$, $\mathsf{N}=2$, and  $\mathsf{n}=10$. The time step size $\delta t$ is $\sfrac{\Delta}{10}$.}
  \label{tbl:multi-stage-wce-advdiff-noncommut-delta-half-order-n10}
  \begin{tabular}[c]{ccc|cccccc}\hline                  \hline
  $\Delta$  &   $\bar{\varrho}_2^{r,2}(T)$    & order     &  $\bar{\varrho}_2^{r,\infty}(T)$ & order  & CPU time (sec.)  \\ \hline \hline
 1.0e-1  & 4.9310e-05 & --              &2.6723e-05& --              &84.00     \\    \hline
 5.0e-2  & 1.4031e-05 & $\Delta^{1.81}$ &7.3571e-06& $\Delta^{ 1.86}$ &160.50    \\    \hline
 2.0e-2  & 2.9085e-06 & $\Delta^{1.71}$ &1.4171e-06& $\Delta^{ 1.80}$ &391.40    \\    \hline
 1.0e-2  & 9.8015e-07 & $\Delta^{1.57}$ &4.4324e-07& $\Delta^{ 1.68}$ &749.40    \\    \hline
 5.0e-3  & 3.5978e-07 & $\Delta^{1.45}$ &1.5082e-07& $\Delta^{ 1.56}$ &1557.60  \\    \hline
 2.0e-3  & 9.8910e-08 & $\Delta^{1.41}$ &3.8369e-08& $\Delta^{ 1.49}$ &3887.50  \\    \hline      \hline
\end{tabular}
\end{table}

With these three examples, we observe that the convergence order of the recursive multi-stage SCM in the second-order moments
is one for commutative and non-commutative noises.
We verified that
our error estimate for WCE, $\Delta^{\mathsf{N}}+\Delta^2$, is valid for commutative noises, see Examples
\ref{exm:sadv-diff-single} and \ref{exm:sadv-diff-com}; the numerical results
for non-commutative noises, see Example \ref{exm:sadv-diff-noncom}, suggest the errors
are of order $\Delta^{\mathsf{N}}+C\sfrac{\Delta}{\mathsf{n}}$ where $C$ is a constant depending
on the coefficients of the noises.

For stochastic advection-diffusion-reaction equations, different formulations of stochastic products
(Ito-Wick product for WCE, Stratonovich product for SCM) lead to different numerical performance. When
the white noise is in the velocity, the  Ito  formulation  will have stronger diffusion
than  that in the Stratonovich formulation
in the resulting PDE. As stronger diffusion requires more  resolution, the recursive multi-stage WCE  with $\mathsf{N}=1$
may produce less accurate results than those
by the recursive multi-stage  SCM with $\mathsf{L}=2$  with  the same PDE solver under the same resolution,
as shown in the first and the third examples.

To achieve convergence 
of approximations of second moments with first-order  in time step 
$\Delta$, we can use  the recursive multi-stage SCM Algorithm~\ref{algo:sadv-diff-s4-scm-mom} with $\mathsf{L}=2$, $\mathsf{n}=1$ and also the
recursive multi-stage WCE Algorithm~\ref{algo:sadv-diff-s4-mom} with $\mathsf{N}=1$, $\mathsf{n}=1$, as both can outperform each other in
certain cases. For commutative noises, Algorithm~\ref{algo:sadv-diff-s4-mom}  with $\mathsf{N}=2$
is preferable when the number of noises, $q $, is small and hence the number of WCE modes is small so that the computational
cost  would grow slowly. 

We also note that  the errors of Algorithms~\ref{algo:sadv-diff-s4-mom} and~\ref{algo:sadv-diff-s4-scm-mom} depend
on the SPDE coefficients and integration time (cf.   theoretical results of Section~3).
For some SPDEs, the constants at powers of $\Delta$ in the errors can be  very large and, to reach desired levels of accuracy, 
 we need to use very small step size  $\Delta$ or develop numerical algorithms further (e.g.,   higher-order or structure-preserving 
 approximations, see such ideas for  stochastic ordinary differential equations e.g. in \cite{MilTre-B04}).
Further,  in practice,  we need to aim at balancing the three parts (truncation of Wiener processes, functional truncation  of WCE/SCM, and space-time discretizations of the deterministic PDEs appearing in the algorithms) of the  errors of Algorithms~\ref{algo:sadv-diff-s4-mom} and~\ref{algo:sadv-diff-s4-scm-mom} for higher computational efficiency.

 \section{Proofs}  \label{eq:error-estimate-wce-scm-proof}
 \subsection{Proof of Theorem \ref{thm:error-estimate-propagator-local}}
 The idea of the proof is to first establish an  estimate for the one-step ($\Delta=T$) error where the global error can readily derived from.
 We  need the following two lemmas  for the one-step errors. Introduce (cf.  \eqref{eq:sadv-diff-wce-solution})
  \begin{equation}\label{eq:sadv-diff-wce-solution-N}
u_{\mathsf{N}}(t,x)=\sum_{\abs{\alpha}\leq \mathsf{N}, \,\alpha\in\cJ_q }\frac{1}{\sqrt{\alpha!}}\varphi_\alpha(t,x)\xi_\alpha.
\end{equation}
\vskip -10pt
\begin{lem}\label{lm:sadv-diff-propagator-trun-N}
 Let $q=1$ in \eqref{eq:sadv-diff}.  Assume  that $\sigma_{i,1}, a_{i,j}, b_i, c,\nu_1$ belong to $\cC^{r+1}_b(\cD)$ 
and $u_0\in H^{r}(\cD)$, where $r\geq \mathsf{N}+1$.
Let $u$ in  \eqref{eq:sadv-diff-wce-solution} be the solution to \eqref{eq:sadv-diff} and  $u_{\mathsf{N}}$ is from \eqref{eq:sadv-diff-wce-solution-N}.
For $C_1 <\delta_{\cL}$, the following estimate holds
\begin{eqnarray*}
\mean{\norm{u(\Delta,\cdot)-u_{\mathsf{N}}(\Delta,\cdot)}^2}\leq   (C_{\floor{r}}\Delta)^{\mathsf{N}+1} e^{2C_{\cL}\Delta} [\frac{e^{C_{\floor{r}}\Delta}}{(\mathsf{N}+1)!}+ \frac{(C_{\floor{r}}\Delta)^{\floor{r}-\mathsf{N}-1}}{\floor{r}!}
   \frac{\delta_{\cL}}{\delta_{\cL}-C_1}]
\norm{u_0}^2_{H^{\floor{r}}},
\end{eqnarray*}
where the constants $\delta_{\cL}$ and $C_{\cL}$
are from \eqref{eq:determ-strong-parabolic-cond} and  $C_{\floor{r}}$ is from \eqref{eq:constant-adv-diff-estimate}.
\end{lem}

 \begin{lem} \label{lm:sadv-diff-propagator-trun-N-n}
Under the assumptions of Lemma   \ref{lm:sadv-diff-propagator-trun-N}  and   $r\geq \mathsf{N}+2$, we have
\begin{equation*}
\mean{\norm{u_{\mathsf{N},\mathsf{n}}(\Delta,\cdot)-u_{\mathsf{N}}(\Delta,\cdot)}^2}\leq \frac{2\Delta^3}{\mathsf{n}\pi^2}C(\mathsf{N}+2,\cL)  \tilde{C}(\mathsf{N},\cL) C_{\mathsf{N}+2}  e^{2C_{\mathsf{N}+2}\Delta+2C_\cL \Delta}
\norm{u_0}^2_{H^{\mathsf{N}+2}},
\end{equation*}
where $C_{\cL}$ is from \eqref{eq:determ-strong-parabolic-cond}, $C(\mathsf{N}+2,\cL) $ is from \eqref{eq:adv-diff-semigroup-estimate},
$\tilde{C}(\mathsf{N},\cL) $ is from \eqref{eq:adv-diff-estimate-diffusion-norm} and $C_{\mathsf{N}+2}$ is from \eqref{eq:constant-adv-diff-estimate}.
\end{lem}

Using Lemmas \ref{lm:sadv-diff-propagator-trun-N} and \ref{lm:sadv-diff-propagator-trun-N-n}, we can establish
the estimate of the global error  stated in Theorem \ref{thm:error-estimate-propagator-local}.    Specifically, the one-step error is bounded by the sum  of
 $\mean{(u(\Delta)-u_{\mathsf{N}}(\Delta))^2}$ and  $\mean{(u_{\mathsf{N}}(\Delta)-u_{\mathsf{N},\mathsf{n}}(\Delta))^2}$, which are estimated
 in Lemmas  \ref{lm:sadv-diff-propagator-trun-N} and  \ref{lm:sadv-diff-propagator-trun-N-n}. Then,
 the global error is estimated based on the recursion nature of Algorithm \ref{algo:sadv-diff-s4-mom} as in the proof in \cite[Theorem 2.4]{LotMR97},
 which completes the proof of Theorem \ref{thm:error-estimate-propagator-local}.

Now we proceed to proving  Lemmas \ref{lm:sadv-diff-propagator-trun-N} and \ref{lm:sadv-diff-propagator-trun-N-n}.
Let us   denote  by $\mathbf{s}^k$   the ordered set
$(s_1,\cdots,s_k)$ and for  $k\geq1$, denote $d\mathbf{s}^k:=\displaystyle
ds_1\dots \,ds_k,$  and \begin{eqnarray*}
 \int^{(k)}(\cdots)\,d\mathbf{s}^k&=&\int_{0}^\Delta\int_{0}^{s_k}\cdots\int_{0}^{s_2}(\cdots)\,ds_1\dots\,ds_k,\\
\int_{(k)}(\cdots)\,d\mathbf{s}^k&=&\int_0^{\Delta}\int_{s_{1}}^\Delta\cdots \int_{s_{k-1}}^\Delta(\cdots)\,ds_k\dots\,ds_2\,ds_1,
 \end{eqnarray*}
 and $\displaystyle F(\Delta;\mathbf{s}^k;x)=\cT_{\Delta-s_k}\cM \cdots
\cT_{s_2-s_1}\cM \cT_{s_1}u_0(x)$, where $\cM:=\cM_1$.

{\em Proof of Lemma \ref{lm:sadv-diff-propagator-trun-N}}.
It follows from   \eqref{eq:determ-strong-parabolic-cond} and the assumptions on  the coefficients   that
 \eqref{eq:adv-diff-semigroup-estimate} and \eqref{eq:adv-diff-semigroup-int-estimate}
hold, cf. \cite[Section 7.1.3]{Evans-B98}. Also, by the assumption that  $\sigma_{i,1}, \nu_1$ belong to $\cC^{r+1}_b(\cD)$,  it can be
readily checked that \eqref{eq:adv-diff-estimate-advection-norm} holds.

By \eqref{eq:sadv-diff-wce-solution}, \eqref{eq:sadv-diff-wce-solution-N} and orthogonality of $\xi_\alpha$ (see \eqref{eq:random-basis-temporal-white-noise}),  we have
\[\mean{\norm{u(\Delta,\cdot)-u_{\mathsf{N}}(\Delta,\cdot)}^2}=\displaystyle\sum_{k>N}\sum_{\abs{\alpha}=k}\frac{\norm{\varphi_{\alpha}(\Delta,\cdot)}^2}{\alpha!}.\]
\vskip -10pt
 Similar to the proof  of Proposition A.1 in \cite{LotMR97}, we have
\[\sum_{\abs{\alpha}=k}\frac{\varphi_{\alpha}^2(\Delta,x)}{\alpha!}= \displaystyle\int^{(k)}\abs{F(\Delta;\mathbf{s}^k;x)}^2\,d\mathbf{s}^k.\]
Then by the Fubini theorem,
\begin{equation}\label{eq:sadv-diff-wce-solution-N-k}
 \sum_{\abs{\alpha}=k}\frac{\norm{\varphi_{\alpha}(\Delta,\cdot)}^2}{\alpha!}= \int^{(k)}\norm{F(\Delta;\mathbf{s}^k;\cdot)}^2\,d\mathbf{s}^k.
 \end{equation}
 \vskip -5pt

Assume that $r>0$ is a   integer.  When $r>0$ is not an integer, we   use $\floor{r}$ instead.

 Denote
$X_k=\cT_{s_k-s_{k-1}}\cM \cdots\cT_{s_2-s_1}\cM \cT_{s_1}u_0$,
$Y_k=\cM X_{k}$, $k\geq 1$ and also $X=\cT_{\Delta-s_k}Y_k$. Then  $X_k=\cT_{s_k-s_{k-1}}Y_{k-1}$ and
$Y_{k-1}=\cM X_{k-1}$.

By the definition of $F$,  \eqref{eq:adv-diff-semigroup-estimate},  \eqref{eq:adv-diff-estimate-advection-norm} and \eqref{eq:constant-adv-diff-estimate},  we have  for  $r\geq k$:
\begin{eqnarray*}
\norm{F(\Delta;\mathbf{s}^k;\cdot)}^2&\leq& e^{2C_{\cL}(\Delta-s_{k})}\norm{Y_k}^2_{H^0} =e^{2C_{\cL}(\Delta-s_{k})}\norm{\cM_{i_k}X_{k}}^2_{H^0}\\
                       &\leq&  \tilde{C}(0,\cM)e^{2C_{\cL}(\Delta-s_{k})}\norm{X_k}^2_{H^1} \\
                       &\leq&  C_1e^{2C_{\cL}(\Delta-s_{k-1})}\norm{Y_{k-1}}^2_{H^1} \leq \cdots \leq  C^k_k e^{2C_{\cL} \Delta} \norm{u_0}^2_{H^k},
\end{eqnarray*}
where  $\tilde{C}(r-1,\cM)$ is from \eqref{eq:adv-diff-estimate-advection-norm} and $C_k$ is defined in \eqref{eq:constant-adv-diff-estimate}.
We then have
\begin{equation}\label{eq:adv-diff-energy-est-smooth}
\int^{(k)}\norm{F(\Delta;\mathbf{s}^k;\cdot)}^2\,d\mathbf{s}^k\leq C^k_ke^{2C_{\cL}\Delta}\norm{u_0}^2_{H^k}  \int^{(k)}\,d\mathbf{s}^k .
\end{equation}
If $r<k$,  by changing  the integration order and  applying  \eqref{eq:adv-diff-semigroup-estimate}, \eqref{eq:adv-diff-estimate-advection-norm} and  \eqref{eq:adv-diff-semigroup-int-estimate},  we get
\begin{eqnarray*}
&&\int^{(k)}\norm{F(\Delta;\mathbf{s}^k;\cdot)}^2\,d\mathbf{s}^k=\int^{(k)}\norm{X}^2\,d\mathbf{s}^k=\int_{(k)}\norm{X}^2\,d\mathbf{s}^k \\
&\leq&  \int_{(k)} e^{2C_{\cL}(\Delta-s_k)}\norm{Y_k}^2\,d\mathbf{s}^k =  \int_{(k)} e^{2C_{\cL}(\Delta-s_k)}\norm{\cM_{i_k}X_{k}}^2\,d\mathbf{s}^k \\
&\leq&\tilde{C}(0,\cM)\int_{(k)} e^{2C_{\cL}(\Delta-s_k)}\norm{X_{k}}^2_{H^1}\,d\mathbf{s}^k\\
&=&\tilde{C}(0,\cM)\int_{(k-1)} \int_{s_{k-1}}^\Delta e^{2C_{\cL}(\Delta-s_{k})}\norm{X_{k}}^2_{H^1}\,ds_k\,d\mathbf{s}^{k-1}\\
&\leq&  \delta_{\cL}^{-1}C_1 \int_{(k-1)}  e^{2C_{\cL}(\Delta-s_{k-1})}\norm{Y_{k-1}}^2\,d\mathbf{s}^{k-1},
\end{eqnarray*}
where $C_1$ is from  \eqref{eq:constant-adv-diff-estimate}.
 Repeating this procedure and using  \eqref{eq:adv-diff-energy-est-smooth}, we obtain
 \begin{equation}\label{eq:adv-diff-energy-est-nonsmooth}
 \resizebox{.9\hsize}{!}
{$\displaystyle\int^{(k)}\norm{F(t;s^{k};\cdot)}^2\,d\mathbf{s}^k \leq
    \delta_{\cL}^{r-k}C_1^{k-r} \int_{(r)}  e^{2C_{\cL}(\Delta-s_{r})}\norm{Y_{r}}^2\,ds^{r}\leq      \delta_{\cL}^{r-k}C_1^{k-r}C^r_r e^{2C_{\cL}t}\norm{u_0}^2_{H^r}  \int^{(r)}\,d\mathbf{s}^r .
$}\end{equation}

By  \eqref{eq:sadv-diff-wce-solution-N}, \eqref{eq:sadv-diff-wce-solution-N-k}, \eqref{eq:adv-diff-energy-est-smooth} and \eqref{eq:adv-diff-energy-est-nonsmooth},
and $\displaystyle\int^{(k)}\,d\mathbf{s}^k =\frac{\Delta^k}{k!}$, we conclude that, for $r\geq  \mathsf{N}+1$ and $C_1<\delta_{\cL}$,
\begin{eqnarray*}
\mean{\norm{u(\Delta,\cdot)-u_{\mathsf{N}}(\Delta,\cdot)}^2} 
                  &=&\sum_{\mathsf{N}<k\leq r}\int^{(k)}\norm{F(\Delta;\mathbf{s}^k;\cdot)}^2\,d\mathbf{s}^k+\sum_{k>r}\int^{(k)}\norm{F(\Delta;\mathbf{s}^k;\cdot)}^2\,d\mathbf{s}^k\\
                  &\leq & \sum_{\mathsf{N}<k\leq r}  \frac{\Delta^k}{k!} C_r^{k} e^{2C_{\cL}\Delta}\norm{u_0}^2_{H^k}
                           +\frac{\Delta^r}{r!}C_r^{r} e^{2C_{\cL}\Delta}\norm{u_0}^2_{H^r}
                           \sum_{k>r}\delta_{\cL}^{r-k}C_1^{k-r}\\
                  &\leq &(C_r\Delta)^{\mathsf{N}+1}e^{2C_{\cL}\Delta} [ \frac{e^{C_r\Delta}}{(\mathsf{N}+1)!}+ \frac{(C_r \Delta )^{r-\mathsf{N}-1}}{r!}
                           \frac{\delta_{\cL}}{\delta_{\cL}-C_1}]\norm{u_0}^2_{H^r}. \hfill \square
\end{eqnarray*}

\begin{rem}
Lemma \ref{lm:sadv-diff-propagator-trun-N}   holds for $r=\infty$ if
$C_\infty<\infty$.  Based on \eqref{eq:adv-diff-energy-est-smooth},  we can prove     that
\begin{eqnarray*}
\mean{\norm{u(\Delta,\cdot)-u_{\mathsf{N}}(\Delta,\cdot)}^2}\leq\sum_{k\geq \mathsf{N}} \frac{\Delta^k}{k!} C_\infty^ke^{2C_{\cL}\Delta}
\norm{u_0}^2_{H^k}
 \leq (C_\infty\Delta)^{\mathsf{N}+1}e^{2C_{\cL} \Delta}   \frac{e^{C_\infty \Delta }}{(\mathsf{N}+1)!} \norm{u_0}^2_{H^\infty}.
 \end{eqnarray*}
\vskip -5pt
If $r<\infty$, we need to   require  that $C_1<\delta_{\cL}$,  i.e.,  $ \tilde{C}(0,\cM)C(1,\cL)<\delta_{\cL}$. For  example,
$\cL=\bigtriangleup$, $\cM_1=\frac{1}{2}D_1$, for which  $\tilde{C}(0,\cM)C(1,\cL)=\frac{1}{2}<\delta_{\cL}=1$.
\end{rem}

{\em Proof of Lemma \ref{lm:sadv-diff-propagator-trun-N-n}}.
It can be proved as in   \cite[p.447]{LotMR97}  that
\begin{eqnarray}\label{eq:estimate-norm-N-n-represent}
\mean{\abs{u_{\mathsf{N}}(\Delta,\cdot)-u_{\mathsf{N},\mathsf{n}}(\Delta,\cdot)}^2}&=&\sum_{l\geq \mathsf{n}+1}\sum_{k=1}^{\mathsf{N}}\sum_{\abs{\alpha}=k, i_k^{\alpha}=l}\frac{\varphi^2_\alpha(\Delta,\cdot)}{\alpha!},
\end{eqnarray}
where  $i_{\abs{\alpha}}^{\alpha}$ is    the index of    last non-zero  element of $\alpha$ and    the    last summation  in the right hand side can be bounded by, see  \cite[(3.7)]{LotMR97},
\begin{equation*}
\sum_{\abs{\alpha}=k, i_k^{\alpha}=l}\frac{  \varphi^2_\alpha(\Delta,x)}{\alpha!}
\leq \int^{(k-1)} \abs{\sum_{j=1}^k \int_{s_{j-1}}^{s_{j+1}} F_j (\Delta;\mathbf{s}^k;x)M_l(s_j)\,ds_j}^2
\,d\mathbf{s}_j^{k},
\end{equation*}
where $d\mathbf{s}_j^{k}=\,ds_1\cdots\,ds_{j-1}\,ds_{j+1}\cdots \,ds_k$, $s_0:=0,\, s_{k+1}:=\Delta$,
 $M_l(t)=\int_0^t m_l(s)\,ds$  and
 \begin{eqnarray*}
 F_j(\Delta;\mathbf{s}^k;x)=\frac{\partial F(\Delta;\mathbf{s}^k;x)}{\partial {s_j}}&=&\cT_{\Delta-s_k}\cM \cdots \cT_{s_{j+1}-s_j}\cM
\cL\cT_{s_{j}-s_{j-1}}\cdots \cT_{s_1}u_0(x)\\
&&-\cT_{\Delta-s_k}\cM \cdots\cM \cL\cT_{s_{j+1}-s_j}\cdots
T_{s_1}u_0(x)=: F_j^1+F_j^2.
\end{eqnarray*}
Then by the Fubini theorem and the Cauchy-Schwarz inequality, we have
\begin{eqnarray*}
\sum_{\abs{\alpha}=k, i_k^{\alpha}=l}\frac{\norm{\varphi_\alpha(\Delta,\cdot)}^2}{\alpha!}
&\leq&      k \int^{(k-1)}   \sum_{j=1}^k \int_{s_{j-1}}^{s_{j+1}}\norm{F_j (\Delta;\mathbf{s}^k;\cdot)}^2\,ds_j  \int_{s_{j-1}}^{s_{j+1}}M_l^2(s_j)\,ds_j
\,d\mathbf{s}_j^{k}. 
\end{eqnarray*}
We claim  (see its proof below) that
 \begin{eqnarray}  \label{eq:F_j-estimate}
 \norm{F_j (\Delta;\mathbf{s}^k;\cdot)}^2 &\leq& 2 \max_{1\leq  j\leq k}\norm{F_j^1}^2  \leq    2 C^{k}_{k+2}  \tilde{C}(k,\cL)  C(k+2,\cL) e^{2C_\cL \Delta}   \norm{u_0}^2_{H^{k+2}}.
 \end{eqnarray}
Thus,  by \eqref{eq:F_j-estimate} we have
 \begin{equation}
\resizebox{.90\hsize}{!}
{$\displaystyle
\sum_{\abs{\alpha}=k, i_k^{\alpha}=l}\frac{\norm{\varphi_\alpha(\Delta,\cdot)}^2}{\alpha!}
 \leq 2 k \Delta      C^{k}_{k+2} C(k+2,\cL) \tilde{C}(k,\cL) e^{2C_\cL \Delta}    \norm{u_0}^2_{H^{k+2}}     \int_0^\Delta M_l^2(s)\,ds  \int^{(k-1)}    \,d\mathbf{s}_j^{k}. \label{eq:estimate-h2-step2-}
$}
\end{equation}
Then by    \eqref{eq:estimate-norm-N-n-represent}, \eqref{eq:estimate-h2-step2-}  and   $M_l(t)=\frac{\sqrt{2\Delta}}{(l-1)\pi}\sin(\frac{(l-1)\pi}{\Delta}t)$  (by \eqref{eq:basis}),  we obtain that
 \begin{eqnarray*}
\mean{\norm{u_{\mathsf{N}}(\Delta,\cdot)-u_{\mathsf{N},\mathsf{n}}(\Delta,\cdot)}^2}&\leq&\sum_{l\geq \mathsf{n}+1} \frac{\Delta^2}{(l-1)^2\pi^2}e^{2C_\cL \Delta}\sum_{k=1}^{\mathsf{N}}  C^{k}_{k+2} C(k+2,\cL)  \tilde{C}(k,\cL) \norm{u_0}^2_{H^{k+2}}  \frac{2k\Delta^{k}}{(k-1)!} \\
&\leq& \frac{2\Delta^3}{\mathsf{n}\pi^2}  e^{2C_\cL \Delta}   \sum_{k=1}^{\mathsf{N}}   C^{k}_{k+2}  \tilde{C}(k,\cL) C(k+2,\cL) \norm{u_0}^2_{H^{k+2}}  \frac{k\Delta^{k-1}}{(k-1)!} \\
&\leq& \frac{2\Delta^3}{\mathsf{n}\pi^2}C_{\mathsf{N}+2}  C(\mathsf{N}+2,\cL)  \tilde{C}(\mathsf{N},\cL) e^{2C_{\mathsf{N}+2}\Delta+2C_\cL \Delta}
\norm{u_0}^2_{H^{\mathsf{N}+2}}.
\end{eqnarray*}

It remains to prove  \eqref{eq:F_j-estimate}.  Note that it is sufficient to estimate $\norm{F^1_j}$ due to
  the same structure of the two terms in $F_j(\Delta;\mathbf{s}^k;x)$.
 By the assumption that $a_{i,j}$ $b_i$ and $c$ belongs to $\cC^{\mathsf{N}+3}_b(\cD)$, it can be readily checked that  \eqref{eq:adv-diff-estimate-diffusion-norm} holds with $l\leq \mathsf{N}+1$.
Repeatedly using  \eqref{eq:adv-diff-semigroup-estimate} and \eqref{eq:adv-diff-estimate-advection-norm} gives
 \begin{eqnarray*}
 \norm{F_j^1}^2&=&\norm{\cT_{\Delta-s_k}\cM \cdots \cT_{s_{j+1}-s_j}
\cM \cL\cT_{s_{j}-s_{j-1}}\cdots \cT_{s_1}u_0}^2\\
             &\leq&    e^{2C_\cL(\Delta-s_k)}\norm{\cM \cdots \cT_{s_{j+1}-s_j}
\cM \cL\cT_{s_{j}-s_{j-1}}\cdots \cT_{s_1}u_0}^2  \\
&\leq&   \tilde{C}(0,\cM) e^{2C_\cL(\Delta-s_k)}\norm{\cT_{s_{k}-s_{k-1}}\cdots \cT_{s_{j+1}-s_j}
\cM\cL\cT_{s_{j}-s_{j-1}}\cdots \cT_{s_1}u_0}^2_{H^1} \\
&\leq&   C_1 e^{2C_\cL(\Delta-s_{k-1})}\norm{\cM \cdots \cT_{s_{j+1}-s_j}
\cM \cL\cT_{s_{j}-s_{j-1}}\cdots \cT_{s_1}u_0}^2_{H^1}\\
&\leq& \cdots  \leq C_{k-j}^{k-j}e^{2C_\cL(\Delta-s_j)}\norm{\cM \cL\cT_{s_{j}-s_{j-1}}\cdots \cT_{s_1}u_0}^2_{H^{k-j}} \\
&\leq& C_{k-j}^{k-j}  \tilde{C}(k-j,\cM) e^{2C_\cL(\Delta-s_j)}\norm{\cL\cT_{s_{j}-s_{j-1}}\cdots \cT_{s_1}u_0}^2_{H^{k-j+1}}\\
&\leq& C_{k-j}^{k-j}  \tilde{C}(k-j,\cM)   \tilde{C}(k-j+1,\cL) e^{2C_\cL(\Delta-s_j)}\norm{ \cT_{s_{j}-s_{j-1}}\cdots \cT_{s_1}u_0}^2_{H^{k-j+3}}\\
&\leq& C_{k-j+1}^{k-j+1}  \tilde{C}(k-j+1,\cL)  e^{2C_\cL(\Delta-s_j)}\norm{\cT_{s_{j}-s_{j-1}}\cM \cdots \cT_{s_1}u_0}^2_{H^{k-j+3}}.  
 \end{eqnarray*}
where we have used   \eqref{eq:adv-diff-estimate-diffusion-norm}  in the last  but one line and the fact that $C(k-j+1,\cL)\geq1$.
Similarly, we have
 \begin{eqnarray*}
 \norm{ \cT_{s_{j}-s_{j-1}}\cM \cdots \cT_{s_1}u_0}^2_{H^{k-j+3}}
 &\leq& C(k-j+3,\cL)        C_{k+2}^{j-1}  e^{2C_\cL  s_j }     \norm{u_0}^2_{H^{k+2}}.
\end{eqnarray*}
Thus, we arrive at   \eqref{eq:F_j-estimate}.   This end the proof of Lemma   \ref{lm:sadv-diff-propagator-trun-N-n}. \hfill $\square$


\subsection{Proof of Theorem \ref{thm:error-est-wzs-global-local-weak}}
To prove Theorem \ref{thm:error-est-wzs-global-local-weak}, we need a  probabilistic representation of the solution to  \eqref{eq:sadv-diff}.
Let   $(\set{B_k(s)}, 1\leq k\leq d, \cF_s^B) $ be  a system of one-dimensional standard Wiener processes  on
a complete probability space   $(\Omega^1,  \cF^1, Q)$ and    independent of $w(s)$   on the space
$(\Omega\otimes \Omega^1, \cF\otimes \cF^1, P\otimes Q)$.
Consider the following
  backward stochastic differential equation on $(\Omega^1,  \cF^1, Q)$, for  $0\leq s\leq t$,
\begin{eqnarray}\label{eq:backward-sde-x}
\hat{d}\hat{X}_{t,x}(s)&=&b(\hat{X}_{t,x}(s))\,ds +\sum_{r=1}^d \alpha_r(\hat{X}_{t,x}(s)) \,\hat{d}B_r(s),\quad \hat{X}_{t,x}(t)=x,
\end{eqnarray}
 The symbol ``$\hat{d}$'' means backward integral, see e.g. \cite{Kun82,Roz-B90} for treatment of  backward stochastic integrals. The $d\times d$ matrix $\alpha(x)$ is  defined by
 $\alpha(x)\alpha^\intercal(x)=2a(x)$. Here $a(x)$ and $b(x)$ are  from \eqref{eq:sadv-diff-coefficients}.
Consider the following backward   stochastic differential equation on $(\Omega\otimes \Omega^1, \cF\otimes \cF^1, P\otimes Q)$ for $0\leq s\leq t$,
\begin{equation}\label{eq:backward-sde-y}
\resizebox{.90\hsize}{!}
{$
\hat{d}\hat{Y}_{t,1,x}(s) = c(\hat{X}_{t,x}(s))\hat{Y}_{t,1,x}(s)\,ds+ \sum_{r=1}^{q } \nu_r(\hat{X}_{t,x}(s))\hat{Y}_{t,1,x}(s) ]\,\hat{d}w_r, \quad \hat{Y}_{t,1,x}(t)=1.
$}
\end{equation}
Here   $c(x)$ and $\nu_r(x)$ are  from \eqref{eq:sadv-diff-coefficients}.
When $u_0(x) \in \cC_b^{2}(\cD)$ and $ \alpha(x),\, b(x),  \, c(x), \, \nu_r(x) \in \cC^0_b(\cD)$ and $\sigma_{i,r}=0$, the solution to  \eqref{eq:sadv-diff}-\eqref{eq:sadv-diff-coefficients} can be represented  by, see e.g. \cite{Kun82},
 \begin{equation}\label{eq:conditional-fk-zakai}
u(t,x)=\meanq{u_0(\hat{X}_{t,x}(0))\exp(\sum_{r=1}^{q } \int_0^t \nu_r(\hat{X}_{t,x}(s))\,\hat{d}w_r(s)+
\int_0^t \bar{c}(\hat{X}_{t,x}(s))\,ds)},
\end{equation}
where $\bar{c}(x)=c(x)-\frac{1}{2}\sum_{r=1}^{q } \nu_r^2(x)$.

Here we first establish  the one-step error  \eqref{eq:error-est-wzs-global-local-weak-one-step} and then the global error  \eqref{eq:error-est-wzs-global-local-weak-multi-step}
We follow the recipe of  the proofs  in \cite[Theorem 3.1]{HuKX02}   and \cite[Theorem 4.4]{BudKal97} where     $\mathsf{n}=1$ and $\mathsf{K}>1$.

We need the following mean-square convergence rate for the one-step error.
\begin{prop}[Mean-square convergence]\label{prop:error-est-wzs-global-local}
Assume that $\sigma_{i,r}=0$ and that  the initial condition $u_0$ and the coefficients in \eqref{eq:sadv-diff-coefficients} are in  $\cC^2_b(\cD)$. Let $u(t,x)$ be the solution to \eqref{eq:sadv-diff} and $\tilde{u}_{\Delta,\mathsf{n}}(t,x)$ the solution to
\eqref{eq:sadv-diff-trun-BM}.   Then for any $\varepsilon>0$,
\begin{equation}\label{eq:prop-error-est-wzs-global-local}
  \mean{\abs{u(\Delta,x)- \tilde{u}_{\Delta, \mathsf{n}}(\Delta,x)}^2}\leq C\exp(C\Delta) (\Delta^3+\Delta^2)\mathsf{n}^{-1+\varepsilon},
\end{equation}
where  the constant  $C>0$ is independent of $\mathsf{n}$.
 \end{prop}

{\em Proof}.
The solution to \eqref{eq:sadv-diff-trun-BM} using the spectral truncation of Brownian motion  $w^{(\Delta,\mathsf{n})}_r$ from \eqref{eq:multi-elem-spectral-exp}
can be represented by, see e.g.  \cite{BudKal97,HuKX02},
\begin{equation}\label{eq:wz-trun-bm-Fayeman-Kac}
\resizebox{.999\hsize}{!}
{$
\tilde{u}_{\Delta,\mathsf{n}}(\Delta,x)=\meanq{u_0(\hat{X}_{\Delta,x}(0))\exp(\sum_{r=1}^{q } \int_0^\Delta \nu_r(\hat{X}_{\Delta,x}(s))\,d{w}^{(\Delta,\mathsf{n})}_r(s)+
\int_0^\Delta \bar{c}(\hat{X}_{\Delta,x}(s))\,ds)}.$
}
\end{equation}
Using $e^x-e^y=e^{\theta x+(1-\theta)y}(x-y)$, $0\leq \theta \leq 1$,
boundedness of $\bar{c}(x)$ and $u_0(x)$  and  the Cauchy-Schwarz inequality (twice), we  have  for some $C>0$:
\begin{eqnarray}\label{eq:zakai-wzs-estimate-first-step}
&&\mean{\abs{\tilde{u}_{\Delta,\mathsf{n}}(\Delta,x) -u(\Delta,x)}^2} \\
&=&
 \mean{\left(\meanq{u_0(\hat{X}_{\Delta,x}(0))\exp(\int_0^\Delta  \bar{c}(\hat{X}_{\Delta,x}(s))\,ds))
 \exp(\sum_{r=1}^{q } \int_0^\Delta \nu_r(\hat{X}_{\Delta,x}(s)) [ \theta\,d{w^{(\Delta,\mathsf{n})}_r(s)}+ (1-\theta)\,\hat{d}w_r(s)]) \right.\notag\\
&&\times \left.
(\sum_{r=1}^{q } \int_0^\Delta \nu_r(\hat{X}_{\Delta,x}(s))[\,d  w^{(\Delta,\mathsf{n})}_r(s) -\,\hat{d}w_r(s)])}\right)^2}\notag\\
&\leq& C \exp(C\Delta)
 \mean{\left(\meanq{
 \exp(\sum_{r=1}^{q } \int_0^\Delta \nu_r(\hat{X}_{\Delta,x}(s)) [ \theta\,d{w^{(\Delta,\mathsf{n})}_r(s)}+ (1-\theta)\,\hat{d}w_r(s)]) \right.\notag\\
&&\times \left.
\abs{\sum_{r=1}^{q } \int_0^\Delta \nu_r(\hat{X}_{\Delta,x}(s))[\,d  w^{(\Delta,\mathsf{n})}_r(s) -\,\hat{d}w_r(s)]}}\right)^2}\notag\\
&\leq&   C\exp(C\Delta)
\big(\mean{\meanq{\exp(\sum_{r=1}^{q } \int_0^\Delta 4\nu_r(\hat{X}_{\Delta,x}(s)) [ \theta\,d  {w}^{(\Delta,\mathsf{n})}_r(s) + (1-\theta)\,\hat{d}w_r(s)])}}\big)^{1/2}\notag\\
&&\times
\big(\mean{\meanq{(\sum_{r=1}^{q } \int_0^\Delta \nu_r(\hat{X}_{\Delta,x}(s))[\,d  {w}^{(\Delta,\mathsf{n})}_r(s) -\,\hat{d}w_r(s)] )^4}}\big)^{1/2}. \notag
\end{eqnarray}
Recall that $\mean{\cdot}=\meanp{\cdot}$ is  the expectation with respect to $P$ only. 
Hence, we need to  estimate $\text{I}_1=\big( \mean{\meanq{(\sum_{r=1}^{q } \int_0^\Delta \nu_r(\hat{X}_{\Delta,x}(s))[\,d  {w}^{(\Delta,\mathsf{n})}_r(s) -\,\hat{d}w_r(s)] )^4}} \big)^{1/2}$
and
$$\text{I}_2= \big(\mean{(\meanq{\exp(\sum_{r=1}^{q } \int_0^\Delta 4\nu_r(\hat{X}_{\Delta,x}(s)) [ \theta\,d   {w}^{(\Delta,\mathsf{n})}_r(s) + (1-\theta)\,\hat{d}w_r(s)])}}\big)^{1/2}.$$
We first estimate $\text{I}_1$.
Due to the independence of $B_k$ and $w_r$, and according to  \cite{Oga84} and   \eqref{eq:bm-spectral-exp}, we have
\begin{eqnarray*}\int_0^\Delta \nu_r(\hat{X}_{\Delta,x}(s))\,\hat{d}w_r(s) =\int_0^\Delta \nu_r(\hat{X}_{\Delta,x}(s))\circ\,dw_r(s) = \sum_{i=0}^\infty \xi_{r,i}\int_0^\Delta \nu_r(\hat{X}_{\Delta,x}(s)) m_{r,i}(s)\,ds.
\end{eqnarray*}
Thus by the Fubini theorem, \eqref{eq:bm-spectral-exp} and \eqref{eq:multi-elem-spectral-exp},  we can represent $\text{I}_1$ as
\begin{eqnarray*}
\text{I}_1&=& \big(\meanq{ \mean{\abs{\sum_{r=1}^{q } [\int_0^\Delta \nu_r(\hat{X}_{\Delta,x}(s))\,d  {w}^{(\Delta,\mathsf{n})}_r(s) -\int_0^\Delta \nu_r(\hat{X}_{\Delta,x}(s))\circ\,dw_r(s)]}^{4}} } \big) ^{1/2}\\
&=&  \big(\meanq{\mean{\abs{\sum_{r=1}^{q } \sum_{i=\mathsf{n}+1}^\infty\xi_{r,i}\int_0^\Delta \nu_r(\hat{X}_{\Delta,x}(s))m_{r,i}(s)\,ds}^{4}} } )^{1/2} \\
 &\leq& \big( 3 \meanq{ \big (\sum_{r=1}^{q }     \sum_{i=\mathsf{n}+1}^\infty  (\int_0^\Delta \nu_r(\hat{X}_{\Delta,x}(s))m_{r,i}(s)\,ds)^2\big)^2} \big)^{1/2},
\end{eqnarray*}
where we  have used twice the fact that    $\hat{X}_{\Delta,x}$   are independent of $w _r$  and  $w^{(\Delta,\mathsf{n})}_r$.
Then by standard estimates of $L^2$-projection error (cf.
\cite[(5.1.10)]{CanHQZ-B06}),  we have for $0< \varepsilon<1,$
\begin{equation}\label{eq:estimate-Slobodeckij semi-norm}
\sum_{i=\mathsf{n}+1}^\infty  (\int_0^\Delta \nu_r(\hat{X}_{\Delta,x}(s))m_{r,i}(s)\,ds)^2 \leq  C \Delta^{1-\varepsilon} \mathsf{n}^{-1+\varepsilon}   \abs{\nu_{r}(\hat{X}_{\Delta,x}(\cdot)) }^2_{\frac{1-\varepsilon}{2},2,[0,\Delta]},
\end{equation}
where the Slobodeckij semi-norm $\abs{f}_{\theta,p,[0,\Delta]}$ is defined by
$(\int_0^\Delta \int_0^{\Delta} \frac{\abs{f(x)-f(y)}^p}{\abs{x-y}^{p\theta+1}}\,dx\,dy)^{1/p} $ and
the constant $\Delta^{1-\varepsilon}$  appears due to the  length of  domain, see e.g. \cite[Chapter 5.4]{CanHQZ-B06}.
Thus, we obtain
\vskip -10pt
\begin{eqnarray}\label{eq:zakai-wzs-I-estimate-inter}
\text{I}_1
 \leq  C \Delta^{1-\varepsilon} \mathsf{n}^{-1+\varepsilon} \big( \sum_{r=1}^{q }\meanq{\abs{\nu_{r}(\hat{X}_{\Delta,x}(\cdot)) }^4_{\frac{1-\varepsilon}{2},2,[0,\Delta]}})^{1/2},   \, 0< \varepsilon<1.
\end{eqnarray}
\vskip -10pt
\noindent
By \eqref{eq:backward-sde-x} and the Ito formula, we have
\begin{eqnarray*}
\hat{X}_{\Delta,x}(s) - \hat{X}_{\Delta,x}(s_1) = \int_{s_1}^s b(\hat{X}_{\Delta,x}(s_2))\,ds_2+  \sum_{k=1}^p \alpha_k(\hat{X}_{\Delta,x}(s_1))[B_k(s)-B_k(s_1)] +R(s_1,s),
\end{eqnarray*}
where $\meanq{\abs{R(s_1,s)}^{2l}}\leq  C\abs{s_1-s}^{2l}$ ($l\geq1$) when  $b(x)$ and $\alpha_k(x)$ belong to $\cC_b^{2}(\cD)$.
By  the Lipschitz continuity of $\nu_1$,  the definition of  the Slobodeckij semi-norm, it is not difficult to show that
\begin{equation}\label{eq:zakai-wzs-I-estimate-component}
\meanq{\abs{\nu_{r}(\hat{X}_{\Delta,x}(\cdot)) }^4_{\frac{1-\varepsilon}{2},2,[0,\Delta]}} \leq C  (\Delta^{4+2\varepsilon}+\Delta^{2+2\varepsilon}).
\end{equation}
Thus, by \eqref{eq:zakai-wzs-I-estimate-inter} and \eqref{eq:zakai-wzs-I-estimate-component}, we have
\begin{equation}\label{eq:zakai-wzs-I-estimate}
 \text{I}_1 \leq C (\Delta^{3}+\Delta^2)  \mathsf{n}^{-1+\varepsilon}.
\end{equation}

Now we estimate $\text{I}_2$.  Using the following facts (see  e.g. \cite[Lemma 2.5]{HuKX02}),
\begin{eqnarray*}
\mean{\exp(\sum_{r=1}^{q }\int_0^\Delta 4\nu_r(\hat{X}_{\Delta,x}(s))\,dw_r)} &= & \exp(\sum_{r=1}^{q } 8\int_0^\Delta \nu_r^2(\hat{X}_{\Delta,x}(s))\,ds), \\
\mean{\exp(\sum_{r=1}^{q }\int_0^\Delta 4\nu_r(\hat{X}_{\Delta,x}(s)) \,d{w}^{(\Delta,\mathsf{n})}_r(s) }&\leq & 4
  \exp(\sum_{r=1}^{q }8\int_0^\Delta \nu_r^2(\hat{X}_{\Delta,x}(s))\,ds),
  \end{eqnarray*}
  we have  
              $\text{I}_2\leq 4\exp(C\Delta)$.
From here, \eqref{eq:zakai-wzs-I-estimate}  and \eqref{eq:zakai-wzs-estimate-first-step},
we reach \eqref{eq:prop-error-est-wzs-global-local}. \hfill $\square$

Now we are ready to prove  Theorem  \ref{thm:error-est-wzs-global-local-weak}, i.e., the   convergence in the second  moments.
%
 For simplicity of notation, we consider $q=1$ while  the case $q>1$ can be proved similarly. Denote
 \begin{eqnarray*}
 U_{\Delta,\mathsf{n},\mathsf{m},\theta}(t,x,\mathbf{y})
   &=:& u_0(\hat{X}_{t,x}(0))\exp( \sum_{i=1}^{\mathsf{n}}\nu_{1,i}y_{i}+
   \theta\sum_{j=\mathsf{n}+1}^{\mathsf{m}} \nu_{1,j} y_{j}+\int_0^t \bar{c}(\hat{X}_{t,x}(s))\,ds), \quad \mathsf{m}\geq\mathsf{n}.
\end{eqnarray*}
 where $\nu_{1,i}(t,x)=\int_0^t \nu_1(\hat{X}_{t,x}(s)) m_{i}(s) \,ds$ for  $ i\leq \mathsf{m}$ ($\hat{X}_{t,x}(s)$ is  the solution to \eqref{eq:backward-sde-x}) and $\mathbf{y}=(y_1,\ldots,y_{\mathsf{n}},  y_{\mathsf{n+1}},  \ldots, y_{\mathsf{m}})$.
  Let us write $\tilde{u}_{\Delta,\mathsf{n},\mathsf{m},\theta}(t,x, \Xi) = \meanq{ U_{\Delta,\mathsf{n},\mathsf{m},\theta}(t,x,\Xi) }$, where
 $\Xi= (\xi_1,\ldots,\xi_{\mathsf{n}},   \xi_{\mathsf{n+1}},  \ldots,   \xi_{\mathsf{m}})$.
 With this notation, we have
 \[\tilde{u}_{\Delta,\mathsf{m}}(t,x)=\tilde{u}_{\Delta,\mathsf{n},\mathsf{m},1}(t,x, \Xi) ,\quad \tilde{u}_{\Delta,\mathsf{n}}(t,x)=\tilde{u}_{\Delta,\mathsf{n},\mathsf{m},0}(t,x, \Xi).\]

 For $\mathsf{m}>\mathsf{n}$, by the first-order Taylor expansion, we have
  \begin{eqnarray} \label{eq:wzs-weak-taylor}
&&\abs{\mean{\tilde{u}_{\Delta,\mathsf{m} }^2 (\Delta,x)- \tilde{u}_{\Delta,\mathsf{n} }^2 (\Delta,x)}} \notag\\
&= & |2\sum_{i,j=\mathsf{n}+1}^{\mathsf{m}}  \frac{1}{(\delta_{i,j}+1)} \int_0^1 \theta (1-\theta) \mean{ \tilde{u}_{\Delta,\mathsf{n},\mathsf{m},\theta}(\Delta,x,\Xi)
\meanq{U_{\Delta,\mathsf{n},\mathsf{m},\theta}(\Delta,x, \Xi)
 \nu_{1,i}(t,x)  \nu_{1,j}(t,x)}  \xi_i\xi_j }  \,d\theta  \notag  \\
& &+2\sum_{i,j=\mathsf{n}+1}^{\mathsf{m}} \frac{1}{(\delta_{i,j}+1)}   \int_0^1 \theta (1-\theta) \mean{\meanq{U_{\Delta,\mathsf{n},\mathsf{m},\theta}(\Delta,x, \Xi)
 \nu_{1,i}(t,x)}  \meanq{U_{\Delta,\mathsf{n},\mathsf{m},\theta}(t,x, \Xi)  \nu_{1,j}(t,x)}  \xi_i\xi_j }  \,d\theta | \notag \\
 &\leq & 2 \abs{\int_0^1 (1-\theta)\theta \mean{ \tilde{u}_{\Delta,\mathsf{n},\mathsf{m},\theta}(\Delta,x,\Xi)
 \meanq{ U_{\Delta,\mathsf{n},\mathsf{m},\theta}(\Delta,x,\Xi) \left(\sum_{i=\mathsf{n}+1}^{\mathsf{m}}\nu_{1,i}(\Delta,x) \xi_i  \right)^2  }  } \,d\theta }  \notag \\
& & +2\abs{\int_0^1 (1-\theta)\theta\mean{\left(\sum_{i=\mathsf{n}+1}^{\mathsf{m}}  \meanq{U_{\Delta, \mathsf{n},\mathsf{m},\theta}(\Delta,x, \Xi)
 \nu_{1,i}(\Delta,x)}  \xi_i  \right)^2} \,d\theta},
\end{eqnarray}
 where  $\delta_{i,j}=1$ if $i=j$ and 0 otherwise and   we have used the facts that $\xi_i$, $i>\mathsf{n}$,  are  independent of $\tilde{u}_{\mathsf{n} }(t,x )$ and  $\mean{\xi_i}=0$.

By the Cauchy-Schwarz inequality (twice), we have  for the first term in  \eqref{eq:wzs-weak-taylor}:
\begin{eqnarray}\label{eq:wzs-weak-I}
&&2\abs{\int_0^1 (1-\theta) \theta \mean{ \tilde{u}_{\Delta,\mathsf{n},\mathsf{m},\theta}(\Delta,x,\Xi)
 \meanq{ U_{\Delta,\mathsf{n},\mathsf{m},\theta}(\Delta,x,\Xi) \left(\sum_{i=\mathsf{n}+1}^{\mathsf{m}}\nu_{1,i}(\Delta,x) \xi_i  \right)^2  }  } \,d\theta }  \notag\\
  &\leq& C  (\mean{
\meanq{\left(\sum_{i=\mathsf{n}+1}^{\mathsf{m}}   \nu_{1,i}(\Delta,x) \xi_i  \right)^8  }  } )^{1/4} .
\end{eqnarray}
Here we also used  that  $\mean{ \tilde{u}_{\Delta,\mathsf{n},\mathsf{m},\theta}^2(\Delta,x,\Xi) },  \mean{  (\meanq{ U_{\Delta,\mathsf{n},\mathsf{m},\theta}^2(\Delta,x,\Xi)})^2}\leq C$, which can be readily checked in the same way  as in the proof of Proposition \ref{prop:error-est-wzs-global-local}.

By the Taylor  expansion for $ U_{\Delta,\mathsf{n},\mathsf{m},\theta}(\Delta,x, \mathbf{y})$,  we have
 \begin{eqnarray*}
 U_{\Delta,\mathsf{n},\mathsf{m},\theta}(\Delta,x, \mathbf{y})
 &=&U_{\Delta,\mathsf{n},\mathsf{m},0}(\Delta,x, \mathbf{y}) +    \sum_{i =\mathsf{n}+1}^{\mathsf{m}}   \nu_{1,i}(\Delta,x)     \int_0^1 (1-\theta_1) \theta_1\theta U_{\Delta,\mathsf{n},\mathsf{m},\theta\theta_1}(\Delta,x, \mathbf{y})\,d\theta_1 y_i.
 \end{eqnarray*}
 Then   by the Cauchy-Schwarz inequality (several times) and the fact that $\xi_i$, $i>\mathsf{n}$ are independent of $U_{\Delta,\mathsf{n},\mathsf{m},0}(t,x, \Xi)$, we have   for the second term in  \eqref{eq:wzs-weak-taylor}:
 \begin{eqnarray}\label{eq:wzs-weak-II}
&&2\abs{\int_0^1 (1-\theta)\theta\mean{\left(\sum_{i=\mathsf{n}+1}^{\mathsf{m}}  \meanq{ U_{\Delta,\mathsf{n},\mathsf{m},\theta}(\Delta,x, \Xi)
 \nu_{1,i}(\Delta,x)}  \xi_i  \right)^2} \,d\theta }    \notag\\
  &\leq&
 4\abs{\int_0^1 (1-\theta) \theta\sum_{i=\mathsf{n}+1}^{\mathsf{m}} \mean{\left( \meanq{U_{\Delta,\mathsf{n},\mathsf{m},0}((\Delta,x, \Xi)
 \nu_{1,i}(\Delta,x)}    \right)^2} \,d\theta }   \notag \\
 &&+ 4\abs{\int_0^1 (1-\theta)\theta^3\mean{\left(   \meanq{ \int_0^1 (1-\theta_1) \theta_1 U_{\Delta,\mathsf{n},\mathsf{m},\theta\theta_1}(\Delta,x, \Xi)\,d\theta_1
 (\sum_{i=\mathsf{n}+1}^{\mathsf{m}}\nu_{1,i}(\Delta,x) \xi_i))^2}   \right)^2} \,d\theta }   \notag \\
&\leq& \mean{  \meanq{U^2_{\Delta,\mathsf{n},\mathsf{m},0}((\Delta,x, \Xi)} } \sum_{i=\mathsf{n}+1}^{\mathsf{m}}   \meanq{\nu^2_{1,i}(\Delta,x)}   \notag\\
 &&+ C\abs{\int_0^1 (1-\theta)\theta^3 (\mean{   \meanq{ \int_0^1 (1-\theta_1)\theta_1U^4_{\Delta,\mathsf{n},\mathsf{m},\theta\theta_1}(\Delta,x, \Xi)\,d\theta_1 } )^{1/2}
 \,d\theta }    } (\mean{\meanq{(\sum_{i=\mathsf{n}+1}^{\mathsf{m}}\nu_{1,i}((\Delta,x) \xi_i)^8}})^{1/2}   \notag\\
 &\leq& C \sum_{i=\mathsf{n}+1}^{\mathsf{m}}   \meanq{\nu^2_{1,i}(\Delta,x)}   + C(\mean{\meanq{(\sum_{i=\mathsf{n}+1}^{\mathsf{m}}\nu_{1,i}(\Delta,x) \xi_i)^8}})^{1/2}.
\end{eqnarray}
Here we used   that $ \mean{  \meanq{U^2_{\Delta,\mathsf{n},\mathsf{m},,0}((\Delta,x, \Xi)} },~\mean{   \meanq{   U^4_{\Delta,\mathsf{n},\mathsf{m},\theta\theta_1}(\Delta,x, \Xi)  }  }\leq C$, which can be readily checked in the same way  as in the proof of Proposition \ref{prop:error-est-wzs-global-local}.

 By \eqref{eq:wzs-weak-taylor}, \eqref{eq:wzs-weak-I} and \eqref{eq:wzs-weak-II}, we  have
 \begin{eqnarray}\label{eq:wz-trun-BM-weak}
&&\abs{\mean{\tilde{u}_{\Delta,\mathsf{m}}^2(\Delta,x) -\tilde{u}_{\Delta,\mathsf{n}}^2(\Delta,x) } }\\
&\leq&  C
 \sum_{i=\mathsf{n}+1}^{\mathsf{m}}   \meanq{\nu^2_{1,i}(\Delta,x)}
 + C\mean{\meanq{(\sum_{i=\mathsf{n}+1}^{\mathsf{m}}\nu_{1,i}(\Delta,x)\xi_i)^8}})^{1/4}   + C (\mean{\meanq{(\sum_{i=\mathsf{n}+1}^{\mathsf{m}}\nu_{1,i}(\Delta,x) \xi_i)^8}})^{1/2} \notag.
  \end{eqnarray}
Similar to the proof of \eqref{eq:zakai-wzs-I-estimate-inter}, we have
  \begin{eqnarray}\label{eq:wz-trun-BM-weak-concl-inter}
 \mean{ \meanq{ (\sum_{i=\mathsf{n}+1}^{\mathsf{m}}   \nu_{1,i} (\Delta,x)   \xi_i )^8  } } &\leq& C\meanq{\left(\sum_{i=\mathsf{n}+1}^{\mathsf{m}}   \nu_{1,i} ^2(\Delta,x)   \right)^4  }
 \leq C\Delta^{4(1-\varepsilon)}\mathsf{n}^{4(1-\varepsilon)}\meanq{\abs{\nu_{1}(\hat{X}_{\Delta,x}(\cdot))}^8_{\frac{1-\varepsilon}{2},2,[0,\Delta]}}. \notag
\end{eqnarray}
Similar to the proof of \eqref{eq:zakai-wzs-I-estimate-component},
we can  estimate  $\meanq{\abs{\nu_{1}(\hat{X}_{\Delta,x}(\cdot))}^8_{\frac{1-\varepsilon}{2},2,[0,\Delta]}}$ as follows:
  \begin{equation*}
    \meanq{\abs{\nu_{1}(\hat{X}_{\Delta,x}(\cdot))}_{\frac{1-\varepsilon}{2},2,[0,\Delta]}^8} \leq  C (\Delta^{8+4\varepsilon} + C \Delta^{4+4\varepsilon}),
  \end{equation*}
  and thus
  \begin{equation} \label{eq:wz-half-norm-bounds-1}
\mean{ \meanq{ (\sum_{i=\mathsf{n}+1}^{\mathsf{m}}   \nu_{1,i} (\Delta,x)   \xi_i )^8  } }\leq  C (\Delta^{12} + C \Delta^{8}).
  \end{equation}
  Similarly, we have
    \begin{equation} \label{eq:wz-half-norm-bounds-2}
\mean{ \meanq{ (\sum_{i=\mathsf{n}+1}^{\mathsf{m}}   \nu_{1,i} (\Delta,x)   \xi_i )^2  } }= \sum_{i=\mathsf{n}+1}^{\mathsf{m}}   \meanq{\nu^2_{1,i}(\Delta,x)}  \leq  C (\Delta^{3} + C \Delta^{2}).
  \end{equation}

By  \eqref{eq:wz-trun-BM-weak}, \eqref{eq:wz-half-norm-bounds-1}  and \eqref{eq:wz-half-norm-bounds-2}, we have
  \begin{equation} \label{eq:wz-half-norm-bounds-inter-final}
  \abs{\mean{\tilde{u}_{\Delta,\mathsf{m}}^2(\Delta,x) -\tilde{u}_{\Delta,\mathsf{n}}^2(\Delta,x) } } \leq  C \exp(C\Delta) (\Delta^{6} + \Delta^{2})  \mathsf{n}^{-1+\varepsilon}.
  \end{equation}
  By the triangle inequality and the Cauchy Schwarz inequality, we obtain
  \begin{eqnarray*}
  \abs{\mean{u^2(\Delta,x) - \tilde{u}_{\Delta,\mathsf{n}}^2(\Delta,x) } }  &\leq&    \abs{\mean{u^2(\Delta,x) - \tilde{u}_{\Delta,\mathsf{m}}^2(\Delta,x) } } + \abs{\mean{\tilde{u}_{\Delta,\mathsf{m}}^2(\Delta,x) - \tilde{u}_{\Delta,\mathsf{n}}^2(\Delta,x) } }, \\
  &\leq & C(\mean{  \abs{u(\Delta,x) - \tilde{u}_{\Delta,\mathsf{m}}(\Delta,x) }^2})^{1/2} + \abs{\mean{\tilde{u}_{\Delta,\mathsf{m}}^2(\Delta,x) - \tilde{u}_{\Delta,\mathsf{n}}^2(\Delta,x) } }.
  \end{eqnarray*}
  The   one-step error  \eqref{eq:error-est-wzs-global-local-weak-one-step}  then
follows from  \eqref{eq:wz-half-norm-bounds-inter-final}, Proposition \ref{prop:error-est-wzs-global-local} and  taking $\mathsf{m}$ to $+\infty$.
The global error   \eqref{eq:error-est-wzs-global-local-weak-multi-step}  is estimated from the recursion nature of Algorithm \ref{algo:sadv-diff-s4-scm-mom} as in the proof in \cite[Theorem 2.4]{LotMR97}.   \hfill $\square$


 \subsection{Proof of Theorem \ref{thm:error-est-wzs-global-local-weak-scm}}

 For any $\mathsf{n}$-dimensional function $\varphi(y_{1},\ldots ,y_{\mathsf{n}})$, we denote
\begin{equation*}
I_{\mathsf{n}}\varphi =\frac{1}{(2\pi )^{\mathsf{n}/2}}\int_{\mathbb{R}%
^{\mathsf{n}}}\varphi (y_{1},\ldots ,y_{\mathsf{n}})\exp \left( {-\frac{1}{2}%
\sum_{i=1}^{\mathsf{n}}y_{i}^{2}}\right) \,d\mathbf{y}.
\end{equation*}%
Introduce the integrals
\begin{equation}
I_{1}^{(k)}\varphi =\frac{1}{\sqrt{2\pi }}\int_{\mathbb{R}}\varphi
(y_{1},\ldots ,y_{k},\ldots ,y_{\mathsf{n}})\exp \left( {-\frac{{y_{k}^{2}}}{2}}%
\right) dy_{k},\ \ k=1,\ldots ,\mathsf{n},  \label{Ik1}
\end{equation}%
and their approximations $Q_{n}^{(k)}$ by the corresponding one-dimensional
Gauss-Hermite quadratures with $n$ nodes. Also, let $\mathcal{U}%
_{i_{k}}^{(k)}=Q_{i_{k}}^{(k)}-Q_{i_{k}-1}^{(k)}.$
By the definition of   Smolyak sparse grid and  using the recipe from the proof of
Lemma~3.4 in \cite{NobTW08}, we obtain
\begin{equation}\label{eq:smolyak-int-recursive-component}
I_{\mathsf{n}}\varphi -A(\mathsf{L},\mathsf{n})\varphi =\sum_{l=2}^{\mathsf{n}}S(\mathsf{L},\mathsf{n})\otimes _{k=l+1}^{\mathsf{n}}{I}%
_{1}^{(k)}\varphi +({I}_{1}^{(1)}-Q_{\mathsf{L}}^{(1)})\otimes _{k=2}^{\mathsf{n}}{I}%
_{1}^{(k)}\varphi ,
\end{equation}%
where
\begin{equation}
S(\mathsf{L},l)=\sum_{i_{1}+\cdots +i_{l-1}+i_{l}=\mathsf{L}+l-1}\otimes _{k=1}^{l-1}\mathcal{U%
}_{i_{k}}^{(k)}\otimes ({I}_{1}^{(l)}-Q_{i_{l}}^{(l)}).
\end{equation}

 Denote     by $D^{\alpha}$   the multivariate derivatives  with  respect to $\mathbf{y}$.
According to the proof of Proposition 3.1 in \cite{ZhangTRK13},  we have
\begin{eqnarray}
&&\left\vert S(\mathsf{L},l)\otimes _{n=l+1}^{\mathsf{n}}I_{1}^{(n)}\varphi \right\vert
\label{eq:sg-error-estimates-components-est} \\
&\leq &\sum_{i_{1}+\cdots +i_{l}=\mathsf{L}+l-1}\frac{(3c/2)^{\#G_{l-1}+1}}{(2\pi
)^{(N-\#F_{l-1})/2}}\int_{\mathbb{R}^{\mathsf{n}-\#F_{l-1}}}\left\vert \otimes _{m\in
F_{l-1}}Q_{1}^{(m)}D^{2\alpha _{l}}\varphi (\mathbf{y})\right\vert  \notag \\
&&\times \exp \left( -\sum_{n\in G_{l-1}}\frac{\beta y_{n}^{2}}{2}-\frac{%
\beta y_{l}^{2}}{2}-\sum_{k=l+1}^{\mathsf{n}}\frac{y_{k}^{2}}{2}\right) \prod_{n\in
G_{l-1}}dy_{n}\times \,dy_{l}\ldots dy_{\mathsf{n}},  \notag
\end{eqnarray}%
where the multi-index $\alpha _{l}=(i_{1}-1,\ldots ,i_{l-1}-1,i_{l},0,\ldots
,0)$ with the $m$-th element $\alpha _{l}^{m},$ the sets $%
F_{l-1}=F_{l-1}(\alpha _{l})=\left\{ m:~\alpha _{l}^{m}=0,\text{ }m=1,\ldots
,l-1\right\} $ and $G_{l-1}=G_{l-1}(\alpha _{l})=\left\{ m:~\alpha
_{l}^{m}>0,\text{ }m=1,\ldots ,l-1\right\} $, the symbols $\#F_{l-1}$ and $%
\#G_{l-1}$ stand for the number of elements in the corresponding sets.
 Here   $c>0,~0<\beta <1$ are only related to the Gauss-Hermite quadrature $Q$, and are independent of the
 number of nodes in the    Gauss-Hermite quadrature,  see e.g. \cite[Theorem 2]{MasMon94}.

 {\em Proof of  of Theorem \ref{thm:error-est-wzs-global-local-weak-scm}.}
 Setting $\varphi(y_1, \cdots,y_{\mathsf{n}})=\tilde{u}_{\Delta,\mathsf{n}}^2(t,x,y_1, \cdots,y_{\mathsf{n}})$, we then have that $A(\mathsf{L},\mathsf{n}) \varphi $ is an  approximation of
 the second  moment  of  the
 solution obtained by the sparse grid collocation methods.
Recall from  \eqref{eq:wz-trun-bm-Fayeman-Kac} that $\tilde{u}_{\Delta,\mathsf{n}}(t, x, \mathbf{y})=\meanq{ U_{\Delta,\mathsf{n}}(t,x,\mathbf{y}) }$
where $U_{\Delta,\mathsf{n}}(t,x,\mathbf{y}) =  U_{\Delta,\mathsf{n},\mathsf{m},0}(t,x,\mathbf{y}) $.

Now we estimate  $D^{2\alpha_l}[\tilde{u}_{\Delta,\mathsf{n}}^2(\Delta,x,y_1, \cdots,y_{\mathsf{n}})]$. To this end, we
need  to  first estimate  $D^{\beta_l}[\tilde{u}_{\Delta,\mathsf{n}}(\Delta,x,y_1, \cdots,y_{\mathsf{n}})]$, where $\beta_l\leq 2\alpha_l.$
By \eqref{eq:estimate-Slobodeckij semi-norm}, we have  for $0<\varepsilon<1$,
\[ \nu_{1,k}^2(\Delta,x) \leq  C (\Delta  \max(k-1, 1))^{\varepsilon-1}  \abs{\nu_1(\hat{X}_{\Delta,x}(\cdot))}^2_{\frac{1-\varepsilon}{2},2,[0,\Delta]},\]
we have, by the Cauchy-Schwarz inequality,
\begin{eqnarray}\label{eq:derivatives-wzs-solu-bound}
&&\abs{D^{\beta_l}\tilde{u}_{\Delta,\mathsf{n}}(\Delta,x,\mathbf{y})}  =
\abs{\meanq{U_{\Delta,\mathsf{n}}(\Delta,x,\mathbf{y}) \prod_{k=1}^{l} (\nu_{1,k}(\Delta,x))^{\beta_l^k}} }  \\
&\leq& (\meanq{U^2_{\Delta,\mathsf{n}}(\Delta,x,\mathbf{y})})^{1/2} (\meanq{\prod_{k=1}^{l} (\nu_{1,k}(\Delta,x))^{2\beta_l^k}})^{1/2}   \notag \\
&\leq&  (C\Delta^{1-\varepsilon} )^{\abs{\beta_l}/2} \prod_{k=2}^l (k-1)^{(\varepsilon-1)\beta_l^{k}/2}  (\meanq{U^2_{\Delta,\mathsf{n}}(\Delta,x,\mathbf{y})})^{1/2}
(\meanq{ \abs{\nu_1(\hat{X}_{\Delta,x}(\cdot))}^{ \abs{2\beta_l}}_{\frac{1-\varepsilon}{2},2,[0,\Delta]}})^{1/2}.\notag
\end{eqnarray}

By the chain rule for multivariate functions, we have 
 \begin{eqnarray*}
D^{2\alpha_l}[\tilde{u}_{\Delta,\mathsf{n}}^2(\Delta,x,\mathbf{y})] &=&
\sum_{\beta_l+\gamma_l=2\alpha_l}(2\alpha_l)! \frac{ D^{\beta_l}\tilde{u}_{\Delta,\mathsf{n}}(\Delta,x,\mathbf{y})  }{\beta_l!} \frac{ D^{{\gamma_l}}\tilde{u}_{\Delta,\mathsf{n}}(\Delta,x,\mathbf{y}) }{\gamma_l!},
\end{eqnarray*}
and thus by  \eqref{eq:derivatives-wzs-solu-bound} and  the fact that  $\sum_{\beta_l+\gamma_l=2\alpha_l}\frac{ (2\alpha_l)!   }{\beta_l!\gamma_l!}=2^{2\abs{\alpha_l}-1}$, we have
\begin{eqnarray*}
\abs{D^{2\alpha_l}[\tilde{u}^2_{\Delta,\mathsf{n}}(\Delta,x,\mathbf{y})]} &\leq & 2^{2\abs{\alpha_l}-1}  (C\Delta^{1-\epsilon})^{\abs{\alpha_l}}\meanq{U^2_{\Delta,\mathsf{n}}(\Delta,x,\mathbf{y})}  \prod_{k=2}^l(k-1)^{ ( \varepsilon-1)\alpha_l^{k}} \\
&& \times
\max_{\beta_l+\gamma_l=2\alpha_l}((\meanq{ \abs{\nu_1(\hat{X}_{\Delta,x}(\cdot)) }^{2\abs{\beta_l}}_{\frac{1-\varepsilon}{2},2,[0,\Delta]}})^{1/2} (\meanq{ \abs{\nu_1(\hat{X}_{\Delta,x}(\cdot))}^{2\abs{\gamma_l}}_{\frac{1-\varepsilon}{2},2,[0,\Delta]}})^{1/2}).
\end{eqnarray*}
Similar to  \eqref{eq:wz-half-norm-bounds-1}, we  have
$\meanq{ \abs{\nu_1(\hat{X}_{\Delta,x}(\cdot))}^{2\abs{\beta_l}}_{\frac{1-\varepsilon}{2},2,[0,\Delta]}}\leq C(\Delta^{ \abs{\beta_l}  (2+\varepsilon)} +   \Delta^{\abs{\beta_l}(1+\varepsilon)})$ and 
\begin{equation}\label{eq:derivatives-wzs-solu-bound-1}
\abs{D^{2\alpha_l}[\tilde{u}^2_{\Delta,\mathsf{n}}(\Delta,x, \mathbf{y})]}
  \leq   C(\Delta^{3\abs{\alpha_l}} + \Delta^{2\abs{\alpha_l}}) \prod_{k=2}^l (k-1)^{(\varepsilon-1)\alpha_l^{k}} \meanq{U^2_{\Delta,\mathsf{n}}(\Delta,x,\mathbf{y})}.
 \end{equation}

Then by \eqref{eq:sg-error-estimates-components-est} and \eqref{eq:derivatives-wzs-solu-bound-1},  we
obtain
\begin{eqnarray} \label{eq:eq:sg-error-estimates-components-est-1}
&&\abs{  S(\mathsf{L},l)\otimes _{n=l+1}^{\mathsf{n}}I_{1}^{(n)}\varphi   } \notag\\
&& \leq C(\Delta^{3\mathsf{L}}+\Delta^{2\mathsf{L}}) (1+(3c/2)^{L\wedge
l})\beta ^{-(\mathsf{L}\wedge l)/2}  \mean{\meanq{U^2_{\Delta,\mathsf{n}}(\Delta,x,\mathbf{y})} }
 \sum_{i_{1}+\cdots +i_{l}=\mathsf{L}+l-1}  \prod_{k=2}^{l} (k-1)^{(\varepsilon-1)\alpha_l^k}\notag \\
&&\leq  C(\Delta^{3\mathsf{L}}+\Delta^{2\mathsf{L}}) (1+(3c/2)^{L\wedge l})\beta ^{-(\mathsf{L}\wedge l)/2}  \epsilon^{1-\mathsf{L}} (l-1)^{\mathsf{L}\varepsilon-1}
\end{eqnarray}%
with the constant $C>0$ which does not depend on $\mathsf{n}$, $\varepsilon$, $\mathsf{L},$ $c,$ $%
\beta ,$ and $l.$ In the last line   we used the fact that
  $ \mean{\meanq{U^2_{\Delta,\mathsf{n}}(\Delta,x,\mathbf{y})} } $ is bounded and that
\begin{eqnarray*}
&&\sum_{i_{1}+\cdots +i_{l}=\mathsf{L}+l-1}  \prod_{k=2}^{l} (k-1)^{(\varepsilon-1)\alpha_l^k}= (l-1)^{\varepsilon-1}\sum_{i_{1}+\cdots +i_{l}=\mathsf{L}+l-1}  \prod_{k=2}^{l} (k-1)^{(\varepsilon-1)(i_k-1)}\\
&\leq&  (l-1)^{\varepsilon-1}
(\sum_{k=2}^l (k-1)^{\varepsilon-1})^{\mathsf{L}-1} \leq   (l-1)^{\varepsilon-1} (\varepsilon^{-1}(l-1)^{\varepsilon})^{\mathsf{L}-1}= \varepsilon^{1-\mathsf{L}} (l-1)^{\mathsf{L}\varepsilon-1}.
\end{eqnarray*}

  Then by \eqref{eq:smolyak-int-recursive-component} and \eqref{eq:eq:sg-error-estimates-components-est-1}, we have \vskip -15pt
 \begin{eqnarray*}
 \abs{ I_{\mathsf{n}}\varphi -A(\mathsf{L},\mathsf{n})\varphi} &\leq&
 C(\Delta^{3\mathsf{L}}+\Delta^{2\mathsf{L}}) (1+(3c/2)^{L\wedge \mathsf{n}})\beta ^{-(\mathsf{L}\wedge \mathsf{n})/2}  \epsilon^{1-\mathsf{L}}
 \sum_{l=2}^{\mathsf{n}} (l-1)^{L\varepsilon-1} \\
 &&  +\abs{(I_{1}^{(1)}-Q_{\mathsf{L}}^{(1)})\otimes _{k=2}^{\mathsf{n}}{I}_{1}^{(k)}\varphi}\\
 &\leq& C(\Delta^{3\mathsf{L}}+\Delta^{2\mathsf{L}}) (1+(3c/2)^{\mathsf{L}\wedge \mathsf{n}})\beta ^{-(\mathsf{L}\wedge \mathsf{n})/2}  \varepsilon^{-\mathsf{L}}{\mathsf{L}}^{-1} \mathsf{n}^{\mathsf{L}\varepsilon},
  \end{eqnarray*}
where the   term  in the second line is estimated by  the classical error estimate for the Gauss-Hermite quadrature $Q$, see e.g. \cite{MasMon94},
 and the estimation of derivatives  \eqref{eq:derivatives-wzs-solu-bound-1}.

 The global error is estimated from the recursion nature of Algorithm \ref{algo:sadv-diff-s4-scm-mom} as in the proof in \cite[Theorem 2.4]{LotMR97}.
 \hfill $\square$
 %

 
\section*{Acknowledgement}
The authors would like to thank  the anonymous  referees for their valuable comments.
 MVT was partially supported by the Leverhulme Trust
Fellowship SAF-2012-006 and  by the UK EPSRC grant EP/K031430/1 and is also grateful to ICERM (Brown University,
Providence) for its hospitality. The rest of the authors were supported partially by a OSD/MURI grant FA9550-09-1-0613,
by NSF grant   DMS-1216437 and also
by the Collaboratory on Mathematics for Mesoscopic Modeling of Materials (CM4) which
is sponsored by DOE. BR was also partially supported by ARO grant W911NF-13-1-0012
and NSF  grant DMS-1148284.

\def\polhk#1{\setbox0=\hbox{#1}{\ooalign{\hidewidth
  \lower1.5ex\hbox{`}\hidewidth\crcr\unhbox0}}} \def\cprime{$'$}


\end{document}